\documentclass[10pt]{article}
\usepackage{amssymb}
\usepackage{amsmath}
\usepackage[dvips]{graphics}
\usepackage{epsfig}

\begin{document}
\title{Sampling from a mixture of different groups of coupons}

\author{ \textbf{Aristides V. Doumas$^1$ and Vassilis G. Papanicolaou$^2$}
\\\\
Department of Mathematics
\\
National Technical University of Athens,
\\
Zografou Campus, 157 80, Athens, GREECE
\\\\
$^1${\tt aris.doumas@hotmail.com} \qquad  $^2${\tt papanico@math.ntua.gr}}
\maketitle

\begin{abstract}
A collector samples coupons with replacement from a pool containing $g$ \textit{uniform} groups of coupons, where
``uniform group" means that all coupons in the group are equally likely to occur (while coupons of different groups have different probabilities
to occur).
For each $j = 1, \dots, g$ let $T_j$ be the number of trials needed
to detect Group $j$, namely to collect all $M_j$ coupons belonging to it at least once. We first derive
formulas for the probabilities $P\{T_1 < \cdots < T_g\}$ and
$P\{T_1 = \bigwedge_{j=1}^g T_j\}$. After that, without severe loss of generality, we restrict ourselves to the case $g=2$ and
compute the asymptotics of $P\{T_1 < T_2\}$ as the number of coupons grows to infinity in a certain manner.
Then, we focus on $T := T_1 \vee T_2$, i.e. the number of trials needed to collect all coupons of the pool
(at least once), and determine the asymptotics of $E[T]$ and $V[T]$, as well as the limiting distribution of $T$
(appropriately normalized) as the number of coupons becomes large.
\end{abstract}

\textbf{Keywords.} Coupon collector problems; urn problems; asymptotics, limiting distribution, Gumbel distribution.
\\\\
\textbf{2010 AMS Mathematics Classification.} 60C05; 60F05; 60F99; 60J10.

\section{Introduction of the problem and main results}
Coupon collector problems (CCP's) are a popular class of urn problems due to their mathematical elegance,
as well as their applications in several areas of science, from computer science and biology to linguistics and the social sciences. The original
problem dates back to De Moivre's treatise \textit{De Mensura Sortis} (1712) and Laplace's \textit{Theorie Analytique des Probabilit\'{e}s} (1812). Nevertheless, new variants of CCP keep arising.

In this paper we study the following CCP version: Suppose we sample coupons independently with replacement from a mixture of $g$ groups of coupons.
The first group consists of $M_1$ coupons each of which having probability $p_1$ to occur, the secong group of
$M_2$ coupons each of which
having probability $p_2$ to occur, and so on (all numbers $M_j$, $p_j$, $j = 1, \dots, g$, are assumed strictly
positive). We call ``Group $j$ coupons" the coupons of the $j$-th group. Notice that under our assumptions we must have
\begin{equation}
M_1 p_1 + \cdots + M_g p_g = 1.
\label{APP1}
\end{equation}
Thus, for each $j = 1, \dots, g$ the $j$-th group is a \textit{uniform} family of $M_j$ coupons, where the term ``uniform" indicates that all
coupons of the group have the same probability $p_j$ to occur. For instance, we can visualize Group 1 as a set of
$M_1$ cards of color $1$ (say red), numbered
from $1$ to $M_1$, Group 2 as a set of $M_2$ cards of color $2$ (say green), numbered from $1$ to $M_2$, and so on,
where each card of color $1$ has
probability $p_1$ to occur, each card of color $2$ has probability $p_2$ to occur, and so on.

Suppose we keep drawing coupons one at a time. Naturally, one quantity of interest is the number $T$ of trials (i.e. draws) needed to detect all
$M_1 + \cdots + M_g$ coupons (at least once). Some ``intermediate" quantities having their own interest are
$T_j :=$ the number of trials needed to detect all
Group $j$ coupons, $j = 1, \dots, g$. Clearly, $T$ can be expressed as
\begin{equation}
T = \bigvee_{j=1}^g T_j,
\label{APP0}
\end{equation}
namely the maximum of $T_1, \dots, T_g$.

It is worth mentioning that if we view the coupon sampling process as a sequence $\{C_n\}_{n \geq 1}$ of independent and identically distributed
random variables, where each $C_n$ takes values in $\{1, 2, \dots, (M_1 + \cdots + M_g)\}$, namely the set of all existing coupons,
with $P\{C_n = i_1\} = p_1$ for $i_1 = 1, 2, \dots, M_1$, $P\{C_n = i_2\} = p_2$ for $i_2 = (M_1 + 1), (M_1 + 2), \dots, (M_1 + M_2)$, and so on
(so that, $\{C_n = i\}$ is identified with the event that the type-$i$ coupon is selected at the $n$-th trial),
then $T_j$, $j = 1, \dots, g$, as well as $T$ are stopping times of the ``coupon filtration''
\begin{equation}
\mathcal{F}_n = \sigma(C_1, \dots, C_n),
\qquad
n \geq 1.
\label{APP0b}
\end{equation}

Our first quantities of study are the probabilities $P\{T_1 < \cdots < T_g\}$ and $P\{T_1 = T_{\min}\}$, where
\begin{equation}
T_{\min} := \bigwedge_{j=1}^g T_j,
\label{APP0c}
\end{equation}
namely the minimum of $T_1, \dots, T_g$ (thus $P\{T_1 = T_{\min}\}$ is the
probability that the Group $1$ is the first group to be detected in its entirety).
Notice that the equality $T_j = T_k$ is impossible unless, of course, $j = k$.
%

\smallskip

\textbf{Theorem 1.}
\begin{align}
& P\{T_1 < \cdots < T_g\}=
\nonumber
\\
& K \int_0^{\infty} \cdots \int_0^{t_3}  e^{-(p_g t_g + \cdots + p_2 t_2)} \left(1 - e^{-p_g t_g}\right)^{M_g - 1} \cdots
\left(1 - e^{-p_2 t_2}\right)^{M_2 - 1} \left(1 - e^{-p_1 t_2}\right)^{M_1} \, dt_2 \cdots dt_g,
\label{NAPP3a}
\end{align}
where
\begin{equation}
K = p_2 p_3 \cdots p_g M_2 M_3 \cdots M_g.
\label{NAPP3b}
\end{equation}
Also,
\begin{align}
& P\{T_1 = T_{\min}\}=
\nonumber
\\
& K \int_0^{\infty} \cdots \int_0^{\infty}
e^{-(p_g t_g + \cdots + p_2 t_2)} \left(1 - e^{-p_g t_g}\right)^{M_g - 1} \cdots
\left(1 - e^{-p_2 t_2}\right)^{M_2 - 1} \left[1 - e^{-p_1 (t_2 \wedge \cdots \wedge t_g)}\right]^{M_1} \, dt_2 \cdots dt_g,
\label{NAPP3c}
\end{align}

\smallskip

\textit{Proof}. Following a suggestion of Professor Sheldon M. Ross \cite{Ro1} we prove the formulas by applying the powerful technique of
``Poissonization."

Let $Z(t)$, $t \geq 0$, be a Poisson process with rate $\lambda = 1$. We imagine that each Poisson event associated to this process is a sampled
coupon, so that $Z(t)$ is the number of sampled coupons at time $t$. Next, for $i = 1, \dots, (M_1 + \cdots + M_g)$, let $Z_i(t)$ be the number of
type-$i$ coupons collected at time $t$. Then, the processes $\{Z_i(t)\}_{t \geq 0}$, $i = 1, \dots, (M_1 + \cdots + M_g)$, are independent Poisson processes with rates $p_1$ for $i = 1, \dots, M_1$, $p_2$ for $i = (M_1+1), \dots, (M_1 + M_2)$, ..., and, finally, $p_g$
for $i = (M_1 + \cdots + M_{g-1} + 1), \dots, (M_1 + \cdots + M_g)$ \cite{Ro2}. Of course, $Z(t) = Z_1(t) + \cdots + Z_{M_1 + \cdots + M_g}(t)$.

If $X_i$, $i = 1, \dots, (M_1 + \cdots + M_g)$, denotes the time when the first type-$i$ coupon is collected, i.e. the time of the first
Poisson event of the process $Z_i(t)$, then the variables $X_1, \dots, X_{M_1 + \cdots + M_g}$ are clearly independent (being associated to
independent processes) and exponentially distributed with parameters $p_1$ for $i = 1, \dots, M_1$, $p_2$ for $i = (M_1+1), \dots, (M_1 + M_2)$
and so on. We now set
\begin{equation}
\tilde{T}_1 := \bigvee_{i=1}^{M_1} X_i,
\qquad
\tilde{T}_2 := \bigvee_{i = M_1 + 1}^{M_1 + M_2} X_i,
\qquad \dots, \qquad
\tilde{T}_g := \bigvee_{i = M_1 + \cdots + M_{g-1} + 1}^{M_1 + \cdots + M_g} X_i.
\label{NP1}
\end{equation}
Thus $\tilde{T}_j$, $j = 1, \dots, g$, is the time when all Group $j$ coupons have been detected (at least once) by the process $Z(t)$ and, hence,
\begin{equation}
P\{T_1 < \cdots < T_g\} = P\{\tilde{T}_1 < \cdots < \tilde{T}_g\}
\quad \text{and} \quad
P\{T_1 = T_{\min}\} = P\{\tilde{T}_1 = \tilde{T}_{\min}\},
\label{NP2}
\end{equation}
where, of course, $\tilde{T}_{\min} := \bigwedge_{j=1}^g \tilde{T}_j$.

From the independence of the exponential random variables $X_i$, $i = 1, \dots, (M_1 + \cdots + M_g)$, it follows that the variables
$\tilde{T}_1, \tilde{T}_2, \ldots, \tilde{T}_g$ are also independent and, furthermore, \eqref{NP1} implies
\begin{equation}
F_j(t) := P\left\{\tilde{T}_j \leq t\right\} = \left(1 - e^{-p_j t}\right)^{M_j},
\qquad
t \geq 0,\ \; j = 1, \dots, g,
\label{NP3a}
\end{equation}
and
\begin{equation}
f_j(t) := F_j'(t) = p_j M_j e^{-p_j t} \left(1 - e^{-p_j t}\right)^{M_j - 1},
\qquad
t \geq 0,\ \; j = 1, \dots, g.
\label{NP3b}
\end{equation}
Therefore,
\begin{align}
P\{\tilde{T}_1 < \cdots < \tilde{T}_g\}
& = \int_0^{\infty} \cdots \int_0^{t_3} \int_0^{t_2} f_g(t_g) \cdots f_2(t_2) f_1(t_1) \, dt_1 dt_2 \cdots dt_g
\nonumber
\\
& = \int_0^{\infty} \cdots \int_0^{t_3} f_g(t_g) \cdots f_2(t_2) F_1(t_2) \, dt_2 \cdots dt_g
\nonumber
\end{align}
and, in view of \eqref{NP2}, \eqref{NP3a}, \eqref{NP3b}, and \eqref{NAPP3b} the above formula is equivalent to \eqref{NAPP3a}.
Likewise,
\begin{align}
P\{\tilde{T}_1 = \tilde{T}_{\min}\}
& = \int_0^{\infty} \cdots \int_0^{\infty} \int_0^{t_2 \wedge \cdots \wedge t_g} f_g(t_g) \cdots f_2(t_2) f_1(t_1) \, dt_1 dt_2 \cdots dt_g
\nonumber
\\
& = \int_0^{\infty} \cdots \int_0^{\infty} f_g(t_g) \cdots f_2(t_2) F_1(t_2 \wedge \cdots \wedge t_g) \, dt_2 \cdots dt_g,
\nonumber
\end{align}
which establishes \eqref{NAPP3c}.
\hfill $\blacksquare$

\smallskip

Notice that one consequence of formulas \eqref{NAPP3a}, \eqref{NAPP3b}, and \eqref{NAPP3c} is that the probabilities  $P\{T_1 < \cdots < T_g\}$ and
$P\{T_1 = T_{\min}\}$ depend only on the ratios $p_2/p_1, \ldots, p_g/p_1$.

\smallskip

\textbf{Corollary 1.} For $\ell = 1, \ldots, g$ we have
\begin{align}
& P\{T_{\ell} = T_{\min}\}=
\nonumber
\\
& (-1)^g \sum_{k_g = 1}^{M_g} \cdots \sum_{k_1 = 1}^{M_1} (-1)^{k_1 + \cdots + k_g}
\binom{M_1}{k_1} \cdots \binom{M_g}{k_g}
\frac{k_{\ell} p_{\ell}} {k_1 p_1 + \cdots + k_g p_g}.
\label{APP3}
\end{align}
In particular, for $g=2$ we have
\begin{equation}
P\{T_1 < T_2\} = \sum_{k=1}^{M_2} \sum_{j=1}^{M_1} (-1)^{j+k} \binom{M_1}{j} \binom{M_2}{k} \frac{p_1 j}{p_1 j + p_2 k}
\label{APP2a}
\end{equation}
and
\begin{equation}
P\{T_2 < T_1\} = \sum_{k=1}^{M_2} \sum_{j=1}^{M_1} (-1)^{j+k} \binom{M_1}{j} \binom{M_2}{k} \frac{p_2 k}{p_1 j + p_2 k}.
\label{APP2b}
\end{equation}

\smallskip

\textit{Proof}. It is enough to prove \eqref{APP3} only for the case $\ell = 1$. For $j = 2, \ldots, g$ we have
\begin{equation}
p_j M_j e^{-p_j t_j} \left(1 - e^{-p_j t_j}\right)^{M_j - 1} = -\sum_{k_j=1}^{M_j} (-1)^{k_j} \binom{M_j}{k_j} p_j k_j e^{-k_j p_j t_j},
\label{NC1a}
\end{equation}
while
\begin{equation}
\left[1 - e^{-p_1 (t_2 \wedge \cdots \wedge t_g)}\right]^{M_1}
= \sum_{k_1=0}^{M_1} (-1)^{k_1} \binom{M_1}{k_1} e^{-k_1 p_1 (t_2 \wedge \cdots \wedge t_g)}.
\label{NC1b}
\end{equation}
Substituting \eqref{NC1a}, \eqref{NC1b}, and \eqref{NAPP3b} in \eqref{NAPP3c} yields
\begin{align}
& P\{T_1 = T_{\min}\}=
\nonumber
\\
& (-1)^{g-1} \sum_{k_g=1}^{M_g} \cdots \sum_{k_2=1}^{M_2} \sum_{k_1=0}^{M_1} (-1)^{k_1 + \cdots + k_g}
\binom{M_g}{k_g} \cdots \binom{M_1}{k_1} I,
\label{NC1c}
\end{align}
where
\begin{equation}
I := \int_0^{\infty} \cdots \int_0^{\infty}
k_g p_g \cdots k_2 p_2 \, e^{-(k_g p_g t_g + \cdots + k_2 p_2 t_2)} e^{-k_1 p_1 (t_2 \wedge \cdots \wedge t_g)} dt_2 \cdots dt_g.
\label{NC1d}
\end{equation}
A quick look at \eqref{NC1d} reveals that
\begin{equation}
I = E\left[e^{-k_1 p_1 Y_{\min}}\right],
\label{NC1dd}
\end{equation}
where $Y_{\min}$ is the minimum of the independent exponential random variables $Y_2, \dots, Y_g$ with parameters $(k_2 p_2), \dots, (k_g p_g)$
respectively. Since (as it is well known) $Y_{\min}$ is exponentially distributed with parameter $k_2 p_2 + \cdots + k_g p_g$, it follows from
\eqref{NC1dd} that
\begin{equation}
I = \frac{k_2 p_2 + \cdots + k_g p_g}{k_1 p_1 + (k_2 p_2 + \cdots + k_g p_g)}
\label{NC1e}
\end{equation}
and the substitution of \eqref{NC1e} in \eqref{NC1c} gives
\begin{align}
& P\{T_1 = T_{\min}\}=
\nonumber
\\
& (-1)^{g-1} \sum_{k_g=1}^{M_g} \cdots \sum_{k_2=1}^{M_2} \sum_{k_1=0}^{M_1} (-1)^{k_1 + \cdots + k_g}
\binom{M_g}{k_g} \cdots \binom{M_1}{k_1} \frac{k_2 p_2 + \cdots + k_g p_g}{k_1 p_1 + k_2 p_2 + \cdots + k_g p_g},
\label{NC1f}
\end{align}
or, equivalently,
\begin{align}
& P\{T_1 = T_{\min}\}=
\nonumber
\\
& (-1)^g \sum_{k_g=1}^{M_g} \cdots \sum_{k_2=1}^{M_2} \sum_{k_1=0}^{M_1} (-1)^{k_1 + \cdots + k_g}
\binom{M_g}{k_g} \cdots \binom{M_1}{k_1} \frac{k_1 p_1}{k_1 p_1 + k_2 p_2 + \cdots + k_g p_g}
\nonumber
\\
& -(-1)^g \sum_{k_g=1}^{M_g} \cdots \sum_{k_2=1}^{M_2} \sum_{k_1=0}^{M_1} (-1)^{k_1 + \cdots + k_g}
\binom{M_g}{k_g} \cdots \binom{M_1}{k_1}.
\label{NC1g}
\end{align}
In the first multiple sum of the right-hand side of \eqref{NC1g}, clearly, the value $k_1=0$ of the dummy variable $k_1$ can be omitted since it does not contribute anything to the sum. Hence we may as well take $k_1$ to vary from $1$ to $M_1$ (instead of $0$ to $M_1$). As for the second multiple
sum of the right-hand side of \eqref{NC1g}, just notice that it can be factored as
\begin{equation*}
\left[\sum_{k_g=1}^{M_g} (-1)^{k_g} \binom{M_g}{k_g}\right] \cdots \left[\sum_{k_2=1}^{M_2} (-1)^{k_2} \binom{M_g}{k_g}\right]
\left[\sum_{k_1=0}^{M_1} (-1)^{k_1} \binom{M_1}{k_1}\right],
\end{equation*}
where, obiously, the last factor is equal to $0$ (being the binomial expansion of $(1 - 1)^{M_1}$). Hence the whole multiple sum vanishes, and \eqref{NC1g} reduces to \eqref{APP3} (for $\ell = 1$).

Formulas \eqref{APP2a} and \eqref{APP2b} follow immediately from \eqref{APP3}.
\hfill $\blacksquare$

\smallskip

An alternative derivation of formula \eqref{APP3} is given in the Subsection 5.1 of the Appendix.

\smallskip


We can number the groups so that $p_1 < p_2 < \cdots < p_g$. Then, our CCP problem is stochastically bounded
between two ``extreme" cases where we have only two groups of coupons: (i) one group consisting of $M_1$ coupons each of which having probability $p_1$ to occur and another group consisting of $M_2 + \cdots + M_g$ coupons each of which having probability $p_2$ to occur and (ii) one group consisting of $M_1$ coupons each of which having probability $p_1$ to occur and another group consisting of $M_2 + \cdots + M_g$ coupons each of which having probability $p_g$ to occur. Hence, the case $g=2$ is quite important since it can, at least, provide upper and lower estimates for the more
general case of an arbitrary number of groups. With this in mind, let us spell out an immediate corollary of Theorem 1.

\smallskip

\textbf{Corollary 2.} We have
\begin{equation}
P\{T_1 < T_2\}
= p_2 M_2 \int_0^{\infty} e^{-p_2 t} \left(1 - e^{-p_1 t}\right)^{M_1} \left(1 - e^{-p_2 t}\right)^{M_2 - 1} dt.
\label{APP5}
\end{equation}
Also (by substitution $x = e^{-t}$ in the above integral),
\begin{align}
P\{T_1 < T_2\}
&= p_2 M_2 \int_0^1 \left(1 - x^{p_1}\right)^{M_1} \left(1 - x^{p_2}\right)^{M_2 - 1} x^{p_2 - 1} dx
\nonumber
\\
&= -\int_0^1 \left(1 - x^{p_1}\right)^{M_1} \left[\left(1 - x^{p_2}\right)^{M_2}\right]' dx.
\label{APP4e}
\end{align}

\medskip

The observation that $P\{T_1 < T_2\}$ depends only on the ratio
\begin{equation}
\lambda := \frac{p_2}{p_1}
\label{APP7}
\end{equation}
yields two slightly simplified equivalent versions of \eqref{APP4e},
namely
\begin{equation}
P\{T_1 < T_2\}
= \lambda M_2 \int_0^1 x^{\lambda - 1} \left(1 - x\right)^{M_1} \left(1 - x^{\lambda}\right)^{M_2 - 1} dx
\label{APP4f}
\end{equation}
and
\begin{equation}
P\{T_1 < T_2\}
= M_2 \int_0^1  \left(1 - x^{1/\lambda}\right)^{M_1} \left(1 - x\right)^{M_2 - 1} dx
\label{APP4g}
\end{equation}
(formula \eqref{APP5} too can be simplified a little by using the fact that $P\{T_1 < T_2\}$ depends only on $\lambda$).
%
%

\smallskip

In the sequel, we will assume that $g=2$, namely that we have only two groups of coupons.
Our goal is to understand the behavior of certain quantities as $\lambda = p_2/p_1$ stays fixed, while $M_1$ and $M_2$ become large in such a way that
\begin{equation}
M_1 = \nu_1 M
\qquad \text{and}\qquad
M_2 = \nu_2 M,
\label{APP6}
\end{equation}
where $\nu_1 \geq 1$ and $\nu_2 \geq 1$ are fixed integers, while the integer $M$ is allowed to grow. Notice that, under these assumptions
$p_1$ and $p_2$ depend on $M$. Then, recalling \eqref{APP1}, namely that $M_1 p_1 + M_2 p_2 = 1$, the quantities
\begin{equation}
\alpha_1 := M_1 p_1 = \frac{\nu_1}{\nu_1 + \lambda \nu_2}
\qquad \text{and} \qquad
\alpha_2 := M_2 p_2 = \frac{\lambda \nu_2}{\nu_1 + \lambda \nu_2} = 1 - \alpha_1
\label{APP7a}
\end{equation}
are independent of $M$ too.

In the rest of the paper we study the asymptotic behavior of certain quantities related to $T_1 = T_1(M)$,
$T_2 = T_2(M)$, and $T = T(M) = T_1\vee T_2$, as the integer $M$ grows large. It is notable that our results determine the order of magnitude of
the corresponding quantities for the case of $g$ groups, for any $g > 2$.

In Section 2 we derive the asymptotic formula (Theorem 2)
\begin{equation*}
P\{T_1 < T_2\} \sim \frac{\nu_2 \lambda \Gamma(\lambda)}{\nu_1^{\lambda}} \cdot \frac{1}{M^{\lambda - 1}},
\qquad
M \to \infty,
\end{equation*}
under the assumption that $\lambda$ of \eqref{APP7} is $>1$ (as usual, the notation $f(M) \sim g(M)$ means that $f(M) / g(M) \to 1$ as
$M \to \infty$).

Section 3 contains some auxiliary topics including a key example discussed in Subsection 3.1. These topics are used in Section 4 in order to determine the asymptotic behavior of the expectation, the variance, as well as the distribution of $T_1$, $T_2$, and $T$ as $M \to \infty$. Some indicative
results of Section 4 are:

(i) A formula for the asymptotics of the expectation of $T$
\begin{equation*}
E[T] = (\nu_1 + \lambda \nu_2) M \ln M  + (\nu_1 + \lambda \nu_2)(\gamma + \ln \nu_1) M + O\left(M^{2-\lambda} \ln M \right),
\qquad
M \to \infty.
\end{equation*}
This formula follows immediately from Theorem 6.

(ii) A formula for the asymptotics of the variance of $T$
\begin{equation*}
V\left[T\right] \sim \frac{\pi^2 (\nu_1 + \lambda \nu_2)^2}{6} \, M^2,
\qquad
M \to \infty
\end{equation*}
(this is Corollary 4).

(iii) The limiting distribution of $T$ (appropriately normalized). We have shown that the random variable
\begin{equation*}
\frac{T - (\nu_1 + \lambda \nu_2) M \, \ln M}{(\nu_1 + \lambda \nu_2) M} - \ln \nu_1
\end{equation*}
converges in distribution to the standard Gumbel random variable as $M \to \infty$ (Theorem 8).

The above three results are presented in Subsection 4.2 and hold under the assumption that $\lambda > 1$.

Finally, for the betterment of the flow of the paper, we have added an appendix (Section 5), where we present an alternative proof of Corollary 1
having its own interest, and, also, the derivations of formulas \eqref{AUX30} and \eqref{A6f}.

\section{Asymptotics of $P\{T_1 < T_2\}$}
Equation \eqref{APP4f} can be written as
\begin{equation}
P\{T_1 < T_2\} = \lambda \nu_2 M \, I_M,
\quad \text{where }\;
I_M = \int_0^1 x^{\lambda - 1} \left(1 - x\right)^{\nu_1 M} \left(1 - x^{\lambda}\right)^{\nu_2 M - 1} dx.
\label{B1}
\end{equation}
Thus, the asymptotic behavior of $P\{T_1 < T_2\}$ as $M \to \infty$ reduces to the asymptotic behavior of $I_M$.

For convenience we will assume from now on, without loss of generality, that
\begin{equation}
\lambda = \frac{p_2}{p_1} > 1.
\label{B2}
\end{equation}

Formula \eqref{B1} yields immediately the following upper bound for $I_M$:
\begin{equation}
I_M < \int_0^1 x^{\lambda - 1} \left(1 - x\right)^{\nu_1 M} dx
= B(\lambda, \nu_1 M + 1) = \frac{\Gamma(\lambda) \, \Gamma(\nu_1 M + 1)}{\Gamma(\lambda + \nu_1 M + 1)},
\label{B3}
\end{equation}
where $B( \, \cdot \, ,\cdot)$ and $\Gamma( \, \cdot)$ denote the Beta and Gamma funtion respectively, while an immediate consequence of Stirling's
formula is that
\begin{equation}
\int_0^1 x^{\lambda - 1} \left(1 - x\right)^{\nu_1 M} dx = \frac{\Gamma(\lambda) \, \Gamma(\nu_1 M + 1)}{\Gamma(\lambda + \nu_1 M + 1)}
\sim \frac{\Gamma(\lambda)}{\nu_1^{\lambda}} \cdot \frac{1}{M^{\lambda}},
\qquad
M \to \infty.
\label{B3a}
\end{equation}

Next, we need to find a satisfactory lower bound  for $I_M$. Let $0 < \varepsilon < 1 - (1 / \lambda)$, so that $(1/\lambda) + \varepsilon < 1$.
Then, \eqref{B1} implies
\begin{equation}
I_M > I_M^{\flat} := \int_0^{M^{-(1/\lambda) - \varepsilon}} x^{\lambda - 1} \left(1 - x\right)^{\nu_1 M} \left(1 - x^{\lambda}\right)^{\nu_2 M - 1} dx.
\label{B4}
\end{equation}
For $x \in [0, M^{-(1/\lambda) - \varepsilon}]$, we have
\begin{equation}
0 \geq (\nu_2 M - 1) \ln\left(1 - x^{\lambda}\right) \geq (\nu_2 M - 1) \ln\left(1 - \frac{1}{M^{1 + \lambda \varepsilon}}\right)
= - \frac{\nu_2}{M^{\lambda \varepsilon}} + O\left(\frac{1}{M^{1 + \lambda \varepsilon}}\right)
\label{B5}
\end{equation}
as $M \to \infty$ (uniformly in $x$). Hence
\begin{equation}
1 \geq \left(1 - x^{\lambda}\right)^{\nu_2 M - 1}
\geq \exp\left(- \frac{\nu_2}{M^{\lambda \varepsilon}}\right) \left[1 + O\left(\frac{1}{M^{1 + \lambda \varepsilon}}\right)\right]
= 1 + O\left(\frac{1}{M^{\lambda \varepsilon}}\right).
\label{B6}
\end{equation}
Using \eqref{B6} in \eqref{B4} yields
\begin{equation}
\int_0^{M^{-(1/\lambda) - \varepsilon}} x^{\lambda - 1} \left(1 - x\right)^{\nu_1 M} dx \geq I_M^{\flat}
\geq \left[1 + O\left(\frac{1}{M^{\lambda \varepsilon}}\right)\right]
\int_0^{M^{-(1/\lambda) - \varepsilon}} x^{\lambda - 1} \left(1 - x\right)^{\nu_1 M} dx
\label{B7}
\end{equation}
or
\begin{equation}
I_M^{\flat} \sim \int_0^{M^{-(1/\lambda) - \varepsilon}} x^{\lambda - 1} \left(1 - x\right)^{\nu_1 M} dx,
\qquad
M \to \infty.
\label{B8}
\end{equation}
Finally, we notice that the fact that $(1/\lambda) + \varepsilon < 1$ implies
\begin{equation}
\int_{M^{-(1/\lambda) - \varepsilon}}^1 x^{\lambda - 1} \left(1 - x\right)^{\nu_1 M} dx
< \left(1 - \frac{1}{M^{(1 / \lambda) + \varepsilon}}\right)^{\nu_1 M}
= O\left(\frac{1}{M^r}\right)
\quad \text{for any }\; r > 0,
\label{B9}
\end{equation}
thus, in view of \eqref{B3a}, formulas \eqref{B8} and \eqref{B9} give
\begin{equation}
I_M^{\flat} \sim \int_0^1 x^{\lambda - 1} \left(1 - x\right)^{\nu_1 M} dx
\sim \frac{\Gamma(\lambda)}{\nu_1^{\lambda}} \cdot \frac{1}{M^{\lambda}},
\qquad
M \to \infty.
\label{B10}
\end{equation}
Hence, the combination of \eqref{B3}, \eqref{B3a}, \eqref{B4}, and \eqref{B10} yields
\begin{equation}
I_M \sim \frac{\Gamma(\lambda)}{\nu_1^{\lambda}} \cdot \frac{1}{M^{\lambda}},
\qquad
M \to \infty.
\label{B11}
\end{equation}
Therefore, by applying \eqref{B11} in \eqref{B1} we obtain the following result.

\smallskip

\textbf{Theorem 2.} If $\lambda = p_2 /p_1 > 1$, then
\begin{equation}
P\{T_1 < T_2\} \sim \frac{\nu_2 \lambda \Gamma(\lambda)}{\nu_1^{\lambda}} \cdot \frac{1}{M^{\lambda - 1}}
= \frac{\nu_2 \Gamma(\lambda + 1)}{\nu_1^{\lambda}} \cdot \frac{1}{M^{\lambda - 1}}
= \frac{\nu_2 }{\nu_1}\cdot \frac{\Gamma(\lambda + 1)}{M_1^{\lambda - 1}}
\label{B12}
\end{equation}
as $M \to \infty$.

\smallskip

Notice that, no matter how big the ratio $M_2 / M_1 = \nu_2 / \nu_1$ is, the probability $P\{T_1 < T_2\}$ approaches $0$ as $M \to \infty$, as long as $\lambda$ is bigger than $1$ (even slightly).

\section{Auxiliary material}
Suppose we sample independently with replacement from a pool of $N$ coupons, where the probability of the $j$-th coupon to occur is $q_j$,
$j = 1, \dots, N$ (the $q_n$'s are usually referred as the ``coupon probabilities"). Let $S = S_N$ denote the number of trials needed in order to
detect all $N$ coupons. Obviously, the possible values of $S_N$ are $N, N+1, \dots$
(it is easy to see that $P\{S_N < \infty\} = 1$ as long as $q_j > 0$ for all $j$;
actually, from the generating function $E\left[z^{-S_N}\right]$, as computed in \cite{DPM}, one can easily see that
$P\{S_N = k\}$ deays exponentially as $k \to \infty$).

For the purposes of this paper we will need a formula for the expectation $E[S_N^{(r)}]$ for any real $r > 0$, where
\begin{equation}
s^{(r)} := \frac{\Gamma(s+r)}{\Gamma(s)}
\label{AUX2}
\end{equation}
is the ``natural" extension of the so-called Pochhammer function.

If we denote by $W_j$ the number of trials needed in order to detect the $j$-th coupon,
then, it is clear that $W_j$ is a geometric random variable with parameter $q_j$ and
\begin{equation*}
S_N = \bigvee_{j=1}^N W_j.
\end{equation*}
However, the above formula for $S_N$ is not very useful, since the $W_j$'s are not independent. Instead, we can employ again the
``Poissonization technique" (see, e.g., \cite{Ro2}) in order to get an explicit formulas for $E[S_N^{(r)}]$.

As in the proof of Thorem 1, we take $Z(t)$, $t \geq 0$, to be a Poisson process with rate $\lambda = 1$. We imagine that each Poisson event associated to $Z$ is a collected coupon,
so that $Z(t)$ is the number of detected coupons at time $t$. Next, for $j = 1, \dots, N$, let $Z_j(t)$ be the number of times that the $j$-th
coupon has been detected up to time $t$. Then, the processes $\{Z_j(t)\}_{t \geq 0}$, $j = 1, \dots, N$, are independent Poisson processes with
rates $q_j$
respectively and, of course, $Z(t) = Z_1(t) + \cdots + Z_N(t)$. If $X_j$, $j = 1, \dots, N$, denotes the time of the first
event of the process $Z_j$, then $X_1, \dots, X_N$ are obviously independent (being associated to independent processes), while their maximum
\begin{equation}
X = \bigvee_{j=1}^N X_j
\label{P1}
\end{equation}
is the time when all different coupons have been detected at least once.

Now, for each $j = 1, \dots, N$ the random variable $X_j$ is exponentially distributed with parameter $q_j$, i.e.
\begin{equation}
P\{X_j \leq t\} = 1 - e^{-q_j t},
\qquad
t \geq 0.
\label{P2}
\end{equation}
It follows from (\ref{P1}) and the independence of the $X_j$'s that
\begin{equation}
P\{X \leq t\} = \prod_{j=1}^N \left( 1 - e^{-q_j t} \right),
\qquad
t \geq 0.
\label{P3}
\end{equation}
Next, we observe that $S_N$ and $X$ are related as
\begin{equation}
X = \sum_{k=1}^{S_N} U_k,
\label{P4}
\end{equation}
where $U_1, U_2, \dots$ are the interarrival times of the process $Z$. It is common knowledge that the $U_j$'s are independent and exponentially distributed random variables with parameter $1$. Hence for any integer $m \geq 1$ the
sum $U_1 + \cdots + U_m$ follows the Erlang distribution with parameters $m$ and $1$. Therefore,
\begin{equation}
E\left[\phi\left(\sum_{k=1}^m U_k\right)\right]
= \int_0^{\infty} \phi(\xi) \, \frac{\xi^{m-1}}{(m-1)!} \, e^{-\xi} d\xi,
\label{P5}
\end{equation}
where $\phi(x)$ is any (Lebesgue) measurable function on $(0, \infty)$ for which the integral in \eqref{P5} makes sense (i.e. converges absolutely).
Noticing that $S_N$ is independent of the $U_j$'s, formulas (\ref{P4}) and (\ref{P5}) imply
\begin{equation}
E\left[\phi(X) \, | \, S_N\right] = \int_0^{\infty} \phi(\xi) \, \frac{\xi^{S_N-1}}{(S_N-1)!} \, e^{-\xi} d\xi
\label{P6a}
\end{equation}
and, consequently (by taking expectations)
\begin{equation}
E\left[\phi(X)\right] = E\left[\int_0^{\infty} \phi(\xi) \, \frac{\xi^{S_N-1}}{(S_N-1)!} \, e^{-\xi} d\xi\right].
\label{P6}
\end{equation}
If we take $\phi(x) = x^r$ for a fixed real number $r > 0$, then \eqref{P6} becomes
\begin{equation}
E\left[X^r\right] = E\left[\int_0^{\infty} \frac{\xi^{S_N + r-1}}{(S_N-1)!} \, e^{-\xi} d\xi\right]
= E\left[\frac{\Gamma(S_N + r)}{(S_N-1)!}\right] = E\left[S_N^{(r)}\right].
\label{P7}
\end{equation}
Finally, by using \eqref{P3} in \eqref{P7} we obtain the following result.

\smallskip

\textbf{Lemma 1.} For any real number $r > 0$ we have
\begin{equation}
E\left[S_N^{(r)}\right] = E\left[\frac{\Gamma(S_N + r)}{\Gamma(S_N)}\right]
= r \int_0^{\infty} t^{r-1}  \left[1 - \prod_{j=1}^N \bigg(1 - e^{-q_j t}\bigg)\right] dt.
\label{AUX3}
\end{equation}
In particular for $r=1$ we have
\begin{equation}
E\left[S_N\right] = \int_0^{\infty}  \left[1 - \prod_{j=1}^N \bigg(1 - e^{-q_j t}\bigg)\right] dt,
\label{mv}
\end{equation}
while for $r=2$ we have
\begin{equation}
E\left[S_N^{(2)}\right] = E\left[S_N (S_N + 1)\right]
= 2 \int_0^{\infty} t  \left[1 - \prod_{j=1}^N \bigg(1 - e^{-q_j t}\bigg)\right] dt.
\label{14}
\end{equation}

\smallskip

\textbf{Remark 1.} Let the random variables $\Xi_1, \ldots, \Xi_N$ be independent and exponentially distributed with parameters $q_1, \ldots, q_N$
respectively. If
\begin{equation*}
\Xi_{\max} := \bigvee_{j=1}^N \Xi_j,
\end{equation*}
then formula \eqref{AUX3} tells us that for any real number $r > 0$ we have
\begin{equation*}
E\left[S_N^{(r)}\right] = E\left[\frac{\Gamma(S_N + r)}{\Gamma(S_N)}\right]
= E\left[\,\Xi_{\max}^r\right].
\end{equation*}

\smallskip

Let us also notice that by expanding the product inside the integral in \eqref{AUX3} and integrate the resulting sum term by term we obtain the expression
\begin{align}
E\left[ S_N^{(r)} \right] & =
\Gamma(r+1) \sum_{\substack{ J\subset
\{ 1, \dots, ,N\}  \\ J\neq \emptyset }} \frac{(-1)^{|J|-1}}
{\left(\sum_{j \in J} q_j \right)^r}
\nonumber
\\
& = \Gamma(r+1) \sum_{m=1}^N (-1)^{m-1} \sum_{1 \leq j_1 < \cdots < j_m \leq N}
\frac{1}{\left(q_{j_1} + \cdots + q_{j_m}\right)^r},
\label{14a}
\end{align}
where $\left|J\right|$ denotes the cardinality of $J$.

\smallskip

\textbf{Remark 2.} Let us first observe that since $S_N$ is always a positive integer, the quantity
\begin{equation}
\frac{\Gamma(r + S_N)}{\Gamma(r) \, \Gamma(S_N)}
\label{AUXC1}
\end{equation}
makes sense for every $r \in \mathbb{C}$; actually, it is entire in $r$ (the poles of $\Gamma(r + S_N)$ are cancelled by the zeros of
$\Gamma(r)^{-1}$). Now, let us look at the function
\begin{equation}
H(r) := \frac{1}{\Gamma(r)} \, E\left[S_N^{(r)}\right] = E\left[\frac{\Gamma(r + S_N)}{\Gamma(r) \, \Gamma(S_N)}\right]
= \sum_{k = N}^{\infty} \frac{\Gamma(k + r)}{\Gamma(r) \, (k-1)!} \, P\{S_N = k\}.
\label{AUXC2}
\end{equation}
Since (i) $\Gamma(k + r) / (k-1)! \sim k^r$ as $k \to \infty$ (see, e.g., formula \eqref{AUX34S} below) and
(ii) $P\{S_N = k\}$ decays exponentially in $k$, it follows that the series in \eqref{AUXC2}
converges uniformly (and absolutely) in $r$ on any compact subset of the complex plane $\mathbb{C}$. Therefore,
$H(r)$ is an entire function
and consequently formula \eqref{AUXC2} implies that $E[S_N^{(r)}]$ is meromorphic in $r$ whose poles are located at $-N, -(N+1), \dots \,\,$.
Although the fact that $E[S_N^{(r)}]$ is meromorphic also follows from \eqref{14a}, it is not obvious from this formula that there are no
poles at $-1, -2, \dots, -(N-1)$.

\smallskip

Now let us consider the ``uniform" case, namely the case where all $N$ coupons are equally likely to occur, i.e.
\begin{equation}
q_j = \frac{1}{N}
\qquad \text{for }\;
j = 1, \dots, N.
\label{AUX4}
\end{equation}
In this case formula \eqref{AUX3} becomes
\begin{equation}
E\left[S_N^{(r)}\right] = r \int_0^{\infty} t^{r-1}  \left[1 - \bigg(1 - e^{- t / N}\bigg)^N \right] dt.
\label{AUX5a}
\end{equation}
Substituting $t = Ns$ in the above integral gives
\begin{equation}
E\left[S_N^{(r)}\right] = N^r \int_0^{\infty} r s^{r-1}  \left[1 - (1 - e^{- s})^N \right] ds.
\label{AUX5}
\end{equation}
Next, we integrate by parts and get
\begin{equation}
E\left[S_N^{(r)}\right] = N^{r+1} \int_0^{\infty} s^r (1 - e^{- s})^{N-1} e^{-s} ds.
\label{AUX6}
\end{equation}
Then, we make the substitution $s = \ln N - \ln x$ (so that $x = N e^{-s}$) in the integral of \eqref{AUX6} and obtain
\begin{equation}
E\left[S_N^{(r)}\right] = N^r \ln^r N \int_0^N \left(1 - \frac{x}{N}\right)^{N-1} \left(1 - \frac{\ln x}{\ln N}\right)^r dx.
\label{AUX7}
\end{equation}
Starting from \eqref{AUX7}, it can be shown that, for any given $r > 0$, the asymptotic behavior of $E[S_N^{(r)}]$ as $N \to \infty$ is
\begin{equation}
E\left[S_N^{(r)}\right] = N^r (\ln N)^r \sum_{k=0}^n \binom{r}{k} \frac{(-1)^k}{\ln^k N}
\int_0^{\infty} e^{-x} (\ln x)^k dx
+ o\left(\frac{1}{\ln^n N}\right)
\label{AUX30}
\end{equation}
for every $n = 1, 2, \dots \;$. Here, $\binom{r}{k}$ stands for the generalized binomial coefficient (in the sense that $r$ is not necessarily a
positive integer---see formulas \eqref{AUX23a} and \eqref{AUX23b} in Subsection 5.2 of the Appendix).

Intuitively, it is not hard to see why \eqref{AUX7} implies \eqref{AUX30}. However the complete proof is quite long and for this reason is given in
the Subsection 5.2 of the Appendix.

To further simplify \eqref{AUX30}, let us first notice that if we differentiate $k$ times the Gamma function
$\Gamma(z) = \int_0^{\infty} t^{z-1} e^{-t} dt$ and then set $z=1$, we get
\begin{equation}
\Gamma^{(k)}(1) =
\int_0^{\infty} e^{-x} (\ln x)^k dx.
\label{AUX31}
\end{equation}
Of course, $\Gamma^{(0)}(1) = \Gamma(1) = 1$. As for the derivatives $\Gamma^{(k)}(1)$, $k = 1, 2, \dots$, there are some known expressions and
recursions (see, e.g., \cite{B-M} and the references therein). For instance,
\begin{equation}
\Gamma^{(1)}(1) = -\gamma,
\quad
\Gamma^{(2)}(1) = \frac{\pi^2}{6} + \gamma^2,
\quad
\Gamma^{(3)}(1) = -\left[2 \zeta(3) + \frac{\pi^2}{2}\gamma + \gamma^3\right],
\quad \text{etc.},
\label{AUX31a}
\end{equation}
where $\gamma = 0.5772...$ is the Euler (or Euler-Macheroni) constant and $\zeta(\cdot)$ is the Riemann Zeta function.

Using \eqref{AUX31} we can write \eqref{AUX30} as
\begin{equation}
E\left[S_N^{(r)}\right] = N^r (\ln N)^r \left[\sum_{k=0}^n \binom{r}{k} \frac{(-1)^k \, \Gamma^{(k)}(1)}{\ln^k N}
+ o\left(\frac{1}{\ln^n N}\right)\right],
\quad
N \to \infty,
\label{AUX32a}
\end{equation}
for every $n = 1, 2, \dots \;$. Formula \eqref{AUX32a} can be written equivalently as an asymptotic series (for the definition of the asymptotic
series and the associated usage of the symbol $\sim$ see, e.g., \cite{B-O})
\begin{equation}
E\left[S_N^{(r)}\right] \sim N^r (\ln N)^r \sum_{k=0}^{\infty} \binom{r}{k} \frac{(-1)^k \, \Gamma^{(k)}(1)}{\ln^k N},
\qquad
N \to \infty,
\label{AUX32b}
\end{equation}
for any $r > 0$ (of course, if $r$ is an integer, $\binom{r}{k} = 0$ for $k > r$ and the series becomes a finite sum). In particular, the leading
behavior of $E[S_N^{(r)}]$ is
\begin{equation}
E\left[S_N^{(r)}\right] \sim N^r (\ln N)^r,
\qquad
N \to \infty.
\label{AUX33}
\end{equation}
Let us also mention that in the case where $r$ is a positive integer there are more detailed expressions for $E[S_N^{(r)}]$ (see, e.g., \cite{DPM}
and the references therein).

Finally, from \eqref{AUX33} we can easily obtain an asymptotic formula for $E[S_N^r]$ as $N \to \infty$. For a fixed $r > 0$ Stirling's formula
yields
\begin{equation}
s^{(r)} = \frac{\Gamma(s+r)}{\Gamma(s)} \sim s^r,
\qquad
s \to \infty.
\label{AUX34S}
\end{equation}
Since $S_N \geq N$, formula \eqref{AUX34} implies that, for any $\varepsilon > 0$ there is a $N_0 = N_0(\varepsilon)$ such that
\begin{equation}
(1 - \varepsilon) \, S_N^{(r)} \leq S_N^r \leq (1 + \varepsilon) \, S_N^{(r)}
\qquad \text{for any }\;
N \geq N_0,
\label{AUX34}
\end{equation}
and consequently,
\begin{equation}
(1 - \varepsilon) \, E\left[S_N^{(r)}\right] \leq E\left[S_N^r\right] \leq (1 + \varepsilon) \, E\left[S_N^{(r)}\right]
\qquad \text{for any }\;
N \geq N_0,
\label{AUX35}
\end{equation}
i.e., in view of \eqref{AUX33},
\begin{equation}
E\left[S_N^r\right] \sim E\left[S_N^{(r)}\right] \sim N^r (\ln N)^r,
\qquad
N \to \infty,
\label{AUX36}
\end{equation}
for any $r > 0$.

\subsection{A preliminary example}
Suppose our set of coupons is $\{0, 1, \dots, N\}$ with corresponding probabilities
\begin{equation}
q_0 = \theta
\qquad \text{and} \qquad
q_j =\frac{1 - \theta}{N},
\quad
j = 1, \dots N,
\label{A2}
\end{equation}
where $\theta \in (0, 1)$ is a given number. Let $S(\theta) = S(\theta; N)$ be the number of trials needed until all $N+1$ coupons are detected
(thus $S(\theta; N) = S_{N+1}$ under the previous notation). Then, \eqref{mv} gives (in the sequel, the dependence of $S(\theta)$
on $N$ will be suppressed for typographical convenience)
\begin{equation*}
E\left[ S(\theta) \right] = \int_0^{\infty} \left[1 - \left(1 - e^{-\theta t} \right) \left(1 - e^{-(1 - \theta) t / N} \right)^N \right] dt
\end{equation*}
or
\begin{equation}
E\left[ S(\theta) \right] = J_1(N; \theta) + J_2(N; \theta),
\label{A3b}
\end{equation}
where
\begin{equation}
J_1(N; \theta) := \int_0^{\infty} \left[1 -  \left(1 - e^{-(1 - \theta) t / N} \right)^N \right] dt
\label{A3c}
\end{equation}
and
\begin{equation}
J_2(N; \theta) := \int_0^{\infty}  e^{-\theta t} \left(1 - e^{-(1 - \theta) t / N} \right)^N dt.
\label{A3d}
\end{equation}
The integral $J_1(N; \theta)$ of \eqref{A3c} reminds the expectation of $S_N$ in the case where all $N$ coupons are equally likely to occur. This is very easy to see via the substitution $y = 1 - e^{-(1 - \theta) t / N}$ which yields
\begin{equation}
J_1(N; \theta) = \frac{N}{1 - \theta} \int_0^1 \frac{1-y^N}{1-y} \, dy = \frac{N}{1 - \theta} \int_0^1 \left(\sum_{j=0}^{N-1} y^j\right)dy
= \frac{N}{1 - \theta} \sum_{j=1}^N \frac{1}{j},
\label{A3e}
\end{equation}
or
\begin{equation}
J_1(N; \theta) = \frac{N H_N}{1 - \theta},
\qquad \text{where} \qquad
H_N := \sum_{j=1}^N \frac{1}{j}.
\label{A3f}
\end{equation}
The quantity $H_N$ is called the $N$-th harmonic number and its full asymptotic expansion, as $N \to \infty$, is well known (see, e.g., \cite{B-O}):
\begin{equation}
H_N \sim \ln N + \gamma + \frac{1}{2N} + \sum_{k=1}^{\infty} \frac{B_{2k}}{2k} \cdot \frac{1}{N^{2k}},
\label{A3g}
\end{equation}
where $B_m$ is the $m$-th Bernoulli number defined by the formula
\begin{equation}
\frac{z}{e^z - 1} = \sum_{m=1}^{\infty} \frac{B_m}{m!} z^m.
\label{A3gg}
\end{equation}
Since $z (e^z - 1)^{-1} + z/2$ is an even function of $z$, we have that $B_{2k+1} = 0$ for all $k \geq 1$.

Next, let us bound the integral $J_2(N; \theta)$ of \eqref{A3d}. For any fixed $\rho > 0$ formula \eqref{A3d} implies
\begin{align}
J_2(N; \theta) &\leq \int_0^{\rho N} \left(1 - e^{-(1 - \theta) t / N} \right)^N dt
+ \int_{\rho N}^{\infty} e^{-\theta t} dt
\nonumber
\\
&\leq \rho N \left(1 - e^{-(1 - \theta) \rho} \right)^N + \frac{1}{\theta \rho N} \, e^{-\theta \rho N},
\label{A3h}
\end{align}
Hence, there is an $\varepsilon_1 > 0$ such that for any fixed $\varepsilon \in (0, \varepsilon_1)$ we have
\begin{equation}
J_2(N; \theta) = O\left(e^{-\varepsilon N}\right),
\qquad
N \to \infty
\label{A3i}
\end{equation}
($\varepsilon$ is a symbol we recycle).

Using \eqref{A3f} and \eqref{A3i} in \eqref{A3b} yields
\begin{equation}
E\left[ S(\theta) \right] = \frac{N H_N}{1 - \theta} + O\left(e^{-\varepsilon N}\right),
\qquad
N \to \infty
\label{A3j}
\end{equation}
(let us recall that
in the case where all $N$ coupons are equally likely to occur we have $E[S_N] = N H_N$).

The full asymptotic expansion of $E\left[ S(\theta) \right]$ can be obtained immediately by applying \eqref{A3g} in \eqref{A3j}. In particular,  \begin{equation}
E\left[ S(\theta) \right] = \frac{N \ln N}{1 - \theta} + \frac{\gamma N}{1 - \theta} + O(1),
\qquad
N \to \infty.
\label{A3k}
\end{equation}

In the same way we can get the asymptotics of the second rising moment $E\left[ S(\theta)^{(2)} \right]$ of $S(\theta)$.
By \eqref{14} and \eqref{A2} we get
\begin{equation}
E\left[S(\theta)^{(2)}\right]
= \tilde{J}_1(N; \theta) + \tilde{J}_2(N; \theta),
\label{A4a}
\end{equation}
where
\begin{equation}
\tilde{J}_1(N; \theta) := 2\int_0^{\infty} t \left[1 - \left(1 - e^{-(1 - \theta) t / N} \right)^N\right] dt
\label{A4b}
\end{equation}
and
\begin{equation}
\tilde{J}_2(N; \theta) := 2\int_0^{\infty} t e^{-\theta t} \left(1 - e^{-(1 - \theta) t / N} \right)^N dt.
\label{A4c}
\end{equation}
The approach we used to bound $J_2(N; \theta)$ applies to $\tilde{J}_2(N; \theta)$ as well and it implies that
there is an $\varepsilon_2 > 0$ such that for any fixed $\varepsilon \in (0, \varepsilon_2)$ we have
\begin{equation}
\tilde{J}_2(N; \theta) = O\left(e^{-\varepsilon N}\right),
\qquad
N \to \infty.
\label{A4d}
\end{equation}
To calculate $\tilde{J}_1(N; \theta)$ we substitute $s = (1 - \theta) t$ in the integral of \eqref{A4b} and obtain
\begin{equation}
\tilde{J}_1(N; \theta) = \frac{2}{(1 - \theta)^2} \int_0^{\infty} s \left[1 - \left(1 - e^{- s / N} \right)^N\right] ds
\label{A4e}
\end{equation}
The integral in the right-hand side of \eqref{A4e} equals $E[S_N^{(2)}]$, where $S_N$ is the number of trials needed to collect all
coupons in the uniform case where all $N$ coupons are equally likely to occur. Since $r=2$ is a very special value, we can get more precise results
than the ones coming directly from formula \eqref{AUX32b}. Indeed, it is not hard to show (see, e.g., \cite{DP}) that
\begin{equation}
E\left[S_N^{(2)}\right] = N^2 \left(H_N^2 + \sum_{j=1}^N \frac{1}{j^2}\right).
\label{A4f}
\end{equation}
Therefore, \eqref{A4e} becomes
\begin{equation}
\tilde{J}_1(N; \theta) = \frac{N^2}{(1 - \theta)^2} \left(H_N^2 + \sum_{j=1}^N \frac{1}{j^2}\right).
\label{A4g}
\end{equation}
Using \eqref{A4d} and \eqref{A4g} in \eqref{A4a} we finally get
\begin{equation}
E\left[S(\theta)^{(2)}\right] = \frac{N^2}{(1 - \theta)^2} \left(H_N^2 + \sum_{j=1}^N \frac{1}{j^2}\right) + O\left(e^{-\varepsilon N}\right),
\qquad
N \to \infty,
\label{A4h}
\end{equation}
for sufficiently small $\varepsilon > 0$. The full asymptotic behavior of $\sum_{j=1}^N j^{-2}$ is well known (see, e.g., \cite{B-O})
\begin{equation}
\sum_{j=1}^N \frac{1}{j^2} \sim \frac{\pi^2}{6} - \frac{1}{N} + \frac{1}{2N^2} - \sum_{k=1}^{\infty} \frac{B_{2k}}{N^{2k+1}},
\label{A4i}
\end{equation}
hence we can obtain easily the full asymptotic expansion of $E\left[S(\theta)^{(2)}\right]$ by using \eqref{A3g} and \eqref{A4i} in \eqref{A4h}.

For the variance of $S(\theta)$ we have
\begin{equation}
V\left[S(\theta)\right] = E\left[S(\theta)^{(2)}\right] - E\left[S(\theta)\right] - E\left[S(\theta)\right]^2,
\label{A5a}
\end{equation}
hence applying \eqref{A3j} and \eqref{A4h} in \eqref{A5a} yields
\begin{equation}
V\left[S(\theta)\right]
= \frac{N^2}{(1 - \theta)^2} \sum_{j=1}^N \frac{1}{j^2} - \frac{N H_N}{1 - \theta}
+ O\left(e^{-\varepsilon N}\right),
\qquad
N \to \infty,
\label{A5b}
\end{equation}
for $\varepsilon > 0$ sufficiently small. Again, the full asymptotic expansion of $V[S(\theta)]$ can be obtained immediately with the help
of \eqref{A3g} and \eqref{A4i}. In particular,
\begin{equation}
V\left[S(\theta)\right]
= \frac{\pi^2 N^2}{6 (1 - \theta)^2}\left[1 + O\left(\frac{\ln N}{N} \right)\right],
\qquad
N \to \infty.
\label{A5c}
\end{equation}

In a similar fashion, we can compute the asymptotics of the fractional rising moments of $S(\theta)$. For $r>0$, in view of \eqref{A2}, formula
\eqref{AUX3} becomes
\begin{align}
E\left[S(\theta)^{(r)}\right]
&= r \int_0^{\infty} t^{r-1} \left[1 - \left(1 - e^{-\theta t} \right) \left(1 - e^{-(1 - \theta) t / N} \right)^N \right] dt
\nonumber
\\
&= r \int_0^{\infty} t^{r-1} \left[1 - \left(1 - e^{-(1 - \theta) t / N} \right)^N \right] dt
+ r \int_0^{\infty} t^{r-1} e^{-\theta t} \left(1 - e^{-(1 - \theta) t / N} \right)^N dt,
\label{A4j}
\end{align}
thus, in the same way we got \eqref{A3i}, we can now get
\begin{equation}
E\left[S(\theta)^{(r)}\right]
= r \int_0^{\infty} t^{r-1} \left[1 - \left(1 - e^{-(1 - \theta) t / N} \right)^N \right] dt
+ O\left(e^{-\varepsilon N}\right),
\qquad
N \to \infty,
\label{A7a}
\end{equation}
for $\varepsilon > 0$ sufficiently small. Next, as usual, we substitute $s = (1 - \theta) t$ in the integral of \eqref{A7a} and obtain
\begin{equation}
E\left[S(\theta)^{(r)}\right]
= \frac{1}{(1 - \theta)^r} \int_0^{\infty} r s^{r-1} \left[1 - \left(1 - e^{-s / N} \right)^N \right] ds \,
+ \, O\left(e^{-\varepsilon N}\right),
\quad
N \to \infty,
\label{A7b}
\end{equation}
In view of \eqref{AUX5a}, the integral in the right-hand side of \eqref{A7b} equals $E[S_N^{(r)}]$, where $S_N$ is the number of trials needed to
collect all coupons in the uniform case. Hence, we can use formula \eqref{AUX32b} in \eqref{A7b} and conclude that
\begin{equation}
E\left[S(\theta)^{(r)}\right]
\sim \frac{N^r (\ln N)^r}{(1 - \theta)^r} \sum_{k=0}^{\infty} \binom{r}{k} \frac{(-1)^k \, \Gamma^{(k)}(1)}{\ln^k N},
\qquad
N \to \infty,
\label{A7c}
\end{equation}
for any $r > 0$. In particular,
\begin{equation}
E\left[S(\theta)^{(r)}\right]
\sim \frac{N^r (\ln N)^r}{(1 - \theta)^r},
\qquad
N \to \infty.
\label{A7d}
\end{equation}
Furthermore, since $S(\theta) \geq N+1$, in the same way we obtained \eqref{AUX36}, we can now get
\begin{equation}
E\left[S(\theta)^r\right] \sim \frac{N^r (\ln N)^r}{(1 - \theta)^r},
\qquad
N \to \infty,
\label{A7e}
\end{equation}
for any $r > 0$.

Finally, we will give the limiting distribution of $S(\theta)$ as $N \to \infty$. The formulas for the moments and the variance of $S(\theta)$
suggest that the law of the random variable $(1 - \theta) S(\theta)$ must be very close to the law of $S_N = $ the number of trials needed to detect
all $N$ coupons in the uniform case where all coupons are equally likely to occur.

The limiting distribution of $S_N$ as $N \to \infty$ has been found in 1961 by Erd\H{o}s and R\'{e}nyi \cite{E-R}:
\begin{equation}
\frac{S_N - N \ln N}{N} \, \overset{D}{\longrightarrow } Y
\qquad \text{as }\;
N \to \infty
\label{A6d}
\end{equation}
(the symbol $\overset{D}{\longrightarrow }$ denotes convergence in distribution) where
\begin{equation}
F(y) = P\{Y \leq y\} = \exp\left(-e^{-y}\right),
\qquad
y \in \mathbb{R},
\label{A6e}
\end{equation}
namely $Y$ is a standard Gumbel random variable. Therefore, it is not surprising that
\begin{equation}
\frac{(1 - \theta) S(\theta) - N \ln N}{N} \, \overset{D}{\longrightarrow } Y
\qquad \text{as }\;
N \to \infty,
\label{A6f}
\end{equation}
where, again, $Y$ is a standard Gumbel random variable.

Our proof of formula \eqref{A6f} is based on characteristic functions. The details are given in the Subsection 5.3 of the Appendix.

%
%

\section{The asymptotic behavior of $T_1$, $T_2$, and $T$}

\subsection{The random variables $T_1$ and $T_2$}
If we are only interested in the variable $T_1 = T_1(M)$ alone, namely the number of trials needed to collect all $M_1 = \nu_1 M$ coupons of Group
$1$, then all the coupons of Group $2$ feel the same to us, and consequently we can assume that the Group $2$ consists of only one coupon having
probability $M_2 p_2 = \alpha_2$ to occur (recall \eqref{APP7a}). Under this point of view, the number of trials $S = S_{M_1 + 1}$ needed to detect
the totality of the $M_1 + 1$ existing coupons (i.e. the $M_1$ coupons of Group $1$ plus the single coupon of Group $2$) can be identified with the
variable $S(\theta) = S(\theta; N)$ studied in Subsection 3.1, where $\theta = \alpha_2$ and $N = M_1 = \nu_1 M$. Although in our notation we will
usually suppress the dependence on $M$ for typographical convenience, we should always keep in mind that both $T_1$ and $S$ below depend on the
integer $M$.

Obviously, $T_1 \leq S$ and the event $\{T_1 < S\}$ happens if and only if the Group $2$ coupon occurs last, namely after detecting all $\nu_1 M$
Group $1$ coupons. Therefore,
\begin{equation}
P\{T_1 < S\} \leq P\{\nu_1 M < S\} = (1 - \alpha_2)^{\nu_1 M} = \alpha_1^{\nu_1 M}
\label{D1}
\end{equation}
(the last equality follows from the fact that, in view of \eqref{APP7a}, $\alpha_1 + \alpha_2 = 1$). This is a rather crude estimate of the
probability of $\{T_1 < S\}$, but it will be sufficient for our purpose.

Next, we will estimate the difference $S - T_1$ in the $L_1$ sense. Let us first notice that,
\begin{equation}
S - T_1 = \left(S - T_1\right) \mathbf{1}_{\{T_1 < S\}}.
\label{DD2}
\end{equation}
Then, taking expectations in \eqref{DD2} yields
\begin{equation}
E\left[S - T_1\right] = E\left[\left(S - T_1\right) \mathbf{1}_{\{T_1 < S\}}\right] = E\left[S - T_1 \, | \, T_1 < S\right] P\{T_1 < S\}.
\label{DD3}
\end{equation}
Now, notice that for $ k = 1, 2, \dots$, we have $P\{S - T_1 = k \, | \, T_1 < S\} = \alpha_1^{k-1} \alpha_2$. Thus, the conditional distribution
of $S - T_1$, given $\{T_1 < S\}$, is geometric with parameter $\alpha_2$. Therefore, $E\left[S - T_1 \, | \, T_1 < S\right] = 1 / \alpha_2$ and
\eqref{DD3} becomes
\begin{equation}
E\left[S\right] - E\left[T_1\right]  = E\left[S - T_1\right] = \frac{1}{\alpha_2} \, P\{T_1 < S\}.
\label{DD4}
\end{equation}
which, in view of \eqref{D1}, implies that $S$ and $T_1$ get very close in the $L_1$ sense as $M \to \infty$. As for the asymptotics of $E[T_1]$, we
can use \eqref{A3j} (with $\theta = \alpha_2$ and $N = \nu_1 M$) and \eqref{D1} in \eqref{DD4} and obtain immediately the following result:

\smallskip

\textbf{Theorem 3.} For every sufficiently small $\varepsilon > 0$ we have
\begin{equation}
E\left[T_1\right] = \frac{\nu_1}{\alpha_1} \, M H_{\nu_1 M} \, + \, O\left(e^{-\varepsilon M}\right)
= (\nu_1 + \lambda \nu_2) M H_{\nu_1 M} \, + \, O\left(e^{-\varepsilon M}\right),
\quad
M \to \infty,
\label{D6a}
\end{equation}
where $H_N$ is the $N$-th harmonic number (see \eqref{A3f}).

Likewise,
\begin{equation}
E\left[T_2\right] = \frac{\nu_2}{\alpha_2} \, M H_{\nu_2 M} \, + \, O\left(e^{-\varepsilon M}\right)
= \big(\lambda^{-1} \nu_1 + \nu_2\big) M H_{\nu_2 M} \, + \, O\left(e^{-\varepsilon M}\right),
\quad
M \to \infty.
\label{D6b}
\end{equation}

\smallskip

For example, in view of \eqref{A3g}, formula \eqref{D6a} implies
\begin{equation}
E\left[T_1\right] = (\nu_1 + \lambda \nu_2) M \ln M \, + \, (\nu_1 + \lambda \nu_2) (\gamma + \, \ln \nu_1) M \, + \, \frac{\alpha_1}{2}
\, + \, O\left(\frac{1}{M}\right),
\quad
M \to \infty,
\label{D6c}
\end{equation}
where, recalling \eqref{APP7a}, we have that $\alpha_1 = \nu_1 / (\nu_1 + \lambda \nu_2)$.

We continue by noticing that in a similar way we can also get easily the asymptotics of the second rising moment of $T_1$.
With the help of Schwarz's inequality (and the fact that $S \geq T_1$) we have
\begin{align}
E\left[S^2\right] - E\left[T_1^2\right] &= E\left[S^2 - T_1^2\right] = E\left[\left(S + T_1\right) \left(S - T_1\right)\right]
\nonumber
\\
&\leq E\left[\left(S + T_1\right)^2\right]^{\frac{1}{2}} E\left[\left(S - T_1\right)^2\right]^{\frac{1}{2}}
\leq 2 E\left[S^2\right]^{\frac{1}{2}} E\left[\left(S - T_1\right)^2\right]^{\frac{1}{2}}.
\label{D7a}
\end{align}
Now, \eqref{DD2} implies that
\begin{equation}
\left(S - T_1\right)^2 = \left(S - T_1\right)^2 \mathbf{1}_{\{T_1 < S\}}
\label{D7b}
\end{equation}
and hence, in the spirit of \eqref{DD3} and \eqref{DD4} we can get
\begin{equation}
E\left[\left(S - T_1\right)^2\right] = E\left[\left(S - T_1\right)^2 \, | \, T_1 < S\right] P\{T_1 < S\}
= \frac{1+\alpha_1}{\alpha_2^2} \, P\{T_1 < S\}.
\label{D7c}
\end{equation}
Using \eqref{D7c} in \eqref{D7a} yields
\begin{equation}
E\left[S^2\right] - E\left[T_1^2\right]
\leq \frac{2 \, \sqrt{1+\alpha_1}}{\alpha_2} \, E\left[S^2\right]^{\frac{1}{2}} P\{T_1 < S\}^{\frac{1}{2}}.
\label{D7d}
\end{equation}
Thus, by \eqref{A4h} (with $\theta = \alpha_2$ and $N = \nu_1 M$) and \eqref{D1} we get that the quantity in the left-hand side of \eqref{D7d}
satisfies
\begin{equation}
E\left[S^2\right] - E\left[T_1^2\right] = O\left(e^{-\varepsilon M}\right),
\qquad
M \to \infty,
\label{D7e}
\end{equation}
for $\varepsilon > 0$ sufficiently small.

Therefore, by applying \eqref{A4h} (with $\theta = \alpha_2$ and $N = \nu_1 M$) in \eqref{D7e} together with Theorem 3 and
\eqref{A3j} (with $\theta = \alpha_2$ and $N = \nu_1 M$) we obtain the following result:

\smallskip

\textbf{Theorem 4.} For every sufficiently small $\varepsilon > 0$ we have
\begin{align}
E\left[T_1(T_1 + 1)\right] &= \left(\frac{\nu_1}{\alpha_1}\right)^2  M^2 \left(H_{\nu_1 M}^2 + \sum_{j=1}^{\nu_1 M} \frac{1}{j^2}\right)
+ O\left(e^{-\varepsilon M}\right)
\nonumber
\\
&= (\nu_1 + \lambda \nu_2)^2 M^2 \left(H_{\nu_1 M}^2 + \sum_{j=1}^{\nu_1 M} \frac{1}{j^2}\right) + O\left(e^{-\varepsilon M}\right)
\label{D8a}
\end{align}
as $M \to \infty$. Likewise,
\begin{align}
E\left[T_2(T_2 + 1)\right] &= \left(\frac{\nu_2}{\alpha_2}\right)^2 M^2 \left(H_{\nu_2 M}^2 + \sum_{j=1}^{\nu_2 M} \frac{1}{j^2}\right)
+ O\left(e^{-\varepsilon M}\right)
\nonumber
\\
&= \big(\lambda^{-1} \nu_1 + \nu_2\big)^2 M^2 \left(H_{\nu_2 M}^2 + \sum_{j=1}^{\nu_2 M} \frac{1}{j^2}\right) + O\left(e^{-\varepsilon M}\right)
\label{D8b}
\end{align}
as $M \to \infty$.

\smallskip

From Theorems 3 and 4 we get immediately the following

\smallskip

\textbf{Corollary 3.} For every sufficiently small $\varepsilon > 0$ we have
\begin{align}
V\left[T_1\right] &= \left(\frac{\nu_1}{\alpha_1}\right)^2 \left(\sum_{j=1}^{\nu_1 M} \frac{1}{j^2}\right) M^2
- \frac{\nu_1}{\alpha_1} \, M H_{\nu_1 M} \, + \, O\left(e^{-\varepsilon M}\right)
\nonumber
\\
&= (\nu_1 + \lambda \nu_2)^2 \left(\sum_{j=1}^{\nu_1 M} \frac{1}{j^2}\right) M^2 - (\nu_1 + \lambda \nu_2)M H_{\nu_1 M}
\, + \,O\left(e^{-\varepsilon M}\right)
\label{D9a}
\end{align}
as $M \to \infty$. Likewise,
\begin{align}
V\left[T_2\right] &= \left(\frac{\nu_2}{\alpha_2}\right)^2 \left(\sum_{j=1}^{\nu_2 M} \frac{1}{j^2}\right) M^2
- \frac{\nu_2}{\alpha_2} \, M H_{\nu_2 M} \, + \, O\left(e^{-\varepsilon M}\right)
\nonumber
\\
&= \big(\lambda^{-1} \nu_1 + \nu_2\big)^2 \left(\sum_{j=1}^{\nu_2 M} \frac{1}{j^2}\right) M^2 - \big(\lambda^{-1} \nu_1 + \nu_2\big) M H_{\nu_2 M}
\, + \,O\left(e^{-\varepsilon M}\right)
\label{D9b}
\end{align}
as $M \to \infty$.

\smallskip

In particular,
\begin{equation}
V\left[T_1\right]
= \frac{\pi^2 (\nu_1 + \lambda \nu_2)^2}{6} \, M^2 \left[1 + O\left(\frac{\ln M}{M} \right)\right],
\qquad
M \to \infty.
\label{D9c}
\end{equation}

\smallskip

Let us, also, mention that a similar approach can be use to determine the asymptotics of $E[T_1^{(r)}]$ and $E[T_1^r]$. Indeed,
formulas \eqref{A7d} and \eqref{A7e} (for $\theta = \alpha_2$ and $N = \nu_1 M$, as usual) imply
\begin{equation}
E\left[T_1^r\right] \sim E\left[T_1^{(r)}\right] \sim E\left[S^{(r)}\right]
\sim (\nu_1 + \lambda \nu_2)^r M^r (\ln M)^r,
\quad
M \to \infty,
\qquad
r > 0.
\label{D8c}
\end{equation}
Likewise,
\begin{equation}
E\left[T_2^r\right] \sim E\left[T_2^{(r)}\right] \sim \big(\lambda^{-1} \nu_1 + \nu_2\big)^r M^r (\ln M)^r,
\quad
M \to \infty,
\qquad
r > 0.
\label{D8d}
\end{equation}

Finally, for $\theta = \alpha_2$ and $N = \nu_1 M$ formula \eqref{A6f} becomes
\begin{equation}
\frac{S - (\nu_1 + \lambda \nu_2) M \, (\ln M + \ln \nu_1)}{(\nu_1 + \lambda \nu_2) M} \overset{D}{\longrightarrow } Y
\qquad \text{as }\;
M \to \infty
\label{D10}
\end{equation}
where $Y$ follows the standard Gumbel distribution displayed in \eqref{A6e}. We can rewrite \eqref{D10} as
\begin{equation}
\frac{T_1 - (\nu_1 + \lambda \nu_2) M \, (\ln M + \ln \nu_1)}{(\nu_1 + \lambda \nu_2) M} \, + \, \frac{S - T_1}{(\nu_1 + \lambda \nu_2) M}
\overset{D}{\longrightarrow } Y
\qquad \text{as }\;
M \to \infty.
\label{D10a}
\end{equation}
However, from \eqref{D1} and \eqref{DD4} we have that $S - T_1 \to 0$ in $L_1$ and, therefore in probability (actually it is easy to see by using
\eqref{D1} and \eqref{DD4} and Chebyshev's inequality that, for any $\delta > 0$ we have $\sum_{M=1}^{\infty} P\{ S - T_1 > \delta\} < \infty$, hence $P\{ S - T_1 > \delta \ \text{i.o.}\} = 0$ and the convergence is almost surely). It follows that $S - T_1 \to 0$ in distribution as $M \to \infty$.
Therefore,
\begin{equation}
\frac{S - T_1}{(\nu_1 + \lambda \nu_2) M}
\overset{D}{\longrightarrow } 0
\qquad \text{as }\;
M \to \infty,
\label{D10b}
\end{equation}
hence by combining \eqref{D10a} and \eqref{D10b} we obtain the following theorem regarding the limiting distribution of $T_1$ (and $T_2$):

\smallskip

\textbf{Theorem 5.}
\begin{equation}
\frac{T_1 - (\nu_1 + \lambda \nu_2) M \, \ln M}{(\nu_1 + \lambda \nu_2) M} - \ln \nu_1 \overset{D}{\longrightarrow } Y
\qquad \text{as }\;
M \to \infty
\label{D11a}
\end{equation}
where
\begin{equation}
F(y) = P\{Y \leq y\} = \exp\left(-e^{-y}\right),
\qquad
y \in \mathbb{R},
\label{D11b}
\end{equation}
namely $Y$ is a standard Gumbel random variable.

Likewise,
\begin{equation}
\frac{T_2 - \big(\lambda^{-1} \nu_1 + \nu_2\big) M \, \ln M}{\big(\lambda^{-1} \nu_1 + \nu_2\big) M} - \ln \nu_2 \overset{D}{\longrightarrow } Y
\qquad \text{as }\;
M \to \infty.
\label{D11c}
\end{equation}

\subsection{The random variable $T$}
We are now ready to determine the asymptotic behavior of the variable $T = T_1 \vee T_2$ as $M \to \infty$. Without loss of generality, as in Section 2, we will assume for convenience that
\begin{equation}
\lambda = \frac{p_2}{p_1} > 1.
\label{E1}
\end{equation}
Let us first observe that we can write
\begin{equation}
T - T_1 = T_1 \vee T_2 - T_1 = \left(T_2 - T_1\right) \mathbf{1}_{\{T_1 < T_2\}}.
\label{E2}
\end{equation}
Taking expectations in \eqref{E2} yields
\begin{equation}
E\left[T - T_1\right] = E\left[\left(T_2 - T_1\right) \mathbf{1}_{\{T_1 < T_2\}}\right]
= E\left[T_2 - T_1 \, | \, T_1 < T_2\right] P\{T_1 < T_2\}.
\label{E3}
\end{equation}
From the fact that $T_1$ and $T_2$ are stopping times of the coupon filtration (recall \eqref{APP0b}) we get
\begin{equation}
E\left[T_2 - T_1 \, | \, T_1 < T_2\right] \leq E\left[T_2\right],
\label{E4}
\end{equation}
thus, using \eqref{E4} in \eqref{E3} gives
\begin{equation}
E\left[T\right] - E\left[T_1\right] = E\left[T - T_1\right] \leq E\left[T_2\right] P\{T_1 < T_2\}.
\label{E5}
\end{equation}
Therefore, by invoking Theorems 2 and 3 we obtain

\smallskip

\textbf{Theorem 6.}
\begin{equation}
E[T] = (\nu_1 + \lambda \nu_2) M H_{\nu_1 M} \, + \, O\left(M^{2-\lambda} \ln M \right),
\qquad
M \to \infty,
\label{E6}
\end{equation}
where, as usual, $H_N$ denotes the $N$-th harmonic number.

\smallskip

Since $\lambda > 1$, formula \eqref{E6} together with \eqref{A3g} imply
\begin{equation}
E[T] = (\nu_1 + \lambda \nu_2) M \ln M  + (\nu_1 + \lambda \nu_2)(\gamma + \ln \nu_1) M + O\left(M^{2-\lambda} \ln M \right),
\qquad
M \to \infty.
\label{E7}
\end{equation}
From Theorem 6 we see that the larger the $\lambda$, the more accurate the asymptotic formula for $E[T]$ becomes. The value $\lambda = 2$ is somehow
critical, since if $\lambda > 2$, then \eqref{E6} yields
\begin{equation}
E[T] = (\nu_1 + \lambda \nu_2) M \ln M  + (\nu_1 + \lambda \nu_2)(\gamma + \ln \nu_1) M + \frac{\nu_1 + \lambda \nu_2}{2 \nu_1} + o(1)
\label{E8}
\end{equation}
as $M \to \infty$.

We continue with the asymptotics of the second rising moment of $T$. We will follow the approach used in the previous subsection for
$E[T_1^{(2)}]$. For better estimates, instead of the Schwarz's inequality we use here the more general H\"{o}lder inequality
(and the fact that $T \leq T_1 + T_2$) to get
\begin{align}
E\left[T^2\right] - E\left[T_1^2\right] &= E\left[T^2 - T_1^2\right] = E\left[\left(T + T_1\right) \left(T - T_1\right)\right]
\nonumber
\\
&\leq E\left[\left(2T_1 + T_2\right)^r\right]^{\frac{1}{r}} E\left[\left(T - T_1\right)^s\right]^{\frac{1}{s}},
\label{E9a}
\end{align}
where
\begin{equation}
r > 1
\qquad \text{and} \qquad
s = \frac{r}{r-1}.
\label{E9aa}
\end{equation}
An immediate upper bound of the first factor of the right-hand side of the inequality in \eqref{E9a} is given by the Minkowski inequality:
\begin{equation}
E\left[\left(2T_1 + T_2\right)^r\right]^{\frac{1}{r}} \leq 2 E\left[T_1^r\right]^{\frac{1}{r}} + E\left[T_2^r\right]^{\frac{1}{r}}.
\label{E9b}
\end{equation}
Now, \eqref{E2} implies
\begin{equation}
\left(T - T_1\right)^s = \left(T - T_1\right)^s \mathbf{1}_{\{T_1 < T_2\}}
\label{E9c}
\end{equation}
and hence, in the spirit of \eqref{E4}
\begin{equation}
E\left[\left(T - T_1\right)^s\right] = E\left[\left(T_2 - T_1\right)^s \, | \, T_1 < T_2\right] P\{T_1 < T_2\}
\leq E\left[T_2^s\right] \, P\{T_1 < T_2\}.
\label{E9d}
\end{equation}
Using \eqref{E9b} and \eqref{E9d} in \eqref{E9a} yields
\begin{equation}
E\left[T^2\right] - E\left[T_1^2\right] \leq
\left(2 E\left[T_1^r\right]^{\frac{1}{r}} + E\left[T_2^r\right]^{\frac{1}{r}}\right) E\left[T_2^s\right]^{\frac{1}{s}} P\{T_1 < T_2\}^{\frac{1}{s}}.
\label{E9e}
\end{equation}
Thus, by using \eqref{D8c}, \eqref{D8d}, and the result of Theorem 2 in \eqref{E9e} we obtain
\begin{equation}
E\left[T^2\right] - E\left[T_1^2\right] = O\left(M^{2 - (\lambda - 1) / s} \ln^2 M \right),
\qquad
M \to \infty.
\label{E9f}
\end{equation}
If we had used the Schwarz's inequality, then we would have been forced to take $r = s = 2$. By using H\"{o}lder inequality, we are free to choose
$r$ as large as we wish and, consequently, in view of \eqref{E9aa}, we can take $s$ arbitrarily close to $1$. Thus, formula \eqref{E9f} is valid
for any $s > 1$ and we can write it as
\begin{equation}
E\left[T^2\right] = E\left[T_1^2\right] + O\left(M^{3 - \lambda + \varepsilon}\right),
\qquad
M \to \infty,
\label{E9ff}
\end{equation}
for any $\varepsilon > 0$. Hence, by Theorems 3 and 6 formula \eqref{E9ff} becomes
\begin{equation}
E\left[T^{(2)}\right] = E\left[T_1^{(2)}\right] + O\left(M^{3 - \lambda + \varepsilon}\right),
\qquad
M \to \infty.
\label{E9g}
\end{equation}
Therefore, by using Theorem 4 in \eqref{E9g} we obtain the following result:

\smallskip

\textbf{Theorem 7.} For every $\varepsilon > 0$ we have
\begin{equation}
E\left[T^{(2)}\right] = (\nu_1 + \lambda \nu_2)^2 M^2 \left(H_{\nu_1 M}^2 + \sum_{j=1}^{\nu_1 M} \frac{1}{j^2}\right)
+ O\left(M^{3 - \lambda + \varepsilon}\right),
\qquad
M \to \infty.
\label{E10}
\end{equation}

\smallskip

Notice that, since $\lambda > 1$ and $\varepsilon$ can be taken arbitrarily close to $0$, the exponent $3 - \lambda + \varepsilon$ in the error term
can be always assumed to be less than $2$ (hence formula \eqref{E9g} is meaningful for any $\lambda > 1$). In particular, from \eqref{E10} we can
immediately deduce that
\begin{equation}
E\left[T^2\right] \sim E\left[T^{(2)}\right] \sim (\nu_1 + \lambda \nu_2)^2 M^2 \ln^2 M,
\qquad
M \to \infty
\label{E10a}
\end{equation}
and, furthermore, in a similar manner we can show that for any $r > 0$ we have
\begin{equation}
E\left[T^r\right] \sim E\left[T^{(r)}\right] \sim (\nu_1 + \lambda \nu_2)^r M^r \ln^r M,
\qquad
M \to \infty.
\label{E10b}
\end{equation}

From Theorems 6 and 7 we get the following corollary.

\smallskip

\textbf{Corollary 4.}
\begin{equation}
V\left[T\right] \sim \frac{\pi^2 (\nu_1 + \lambda \nu_2)^2}{6} \, M^2,
\qquad
M \to \infty.
\label{E11}
\end{equation}

\smallskip

Finally, let us determine the limiting distribution of $T$ as $M \to \infty$. Formula \eqref{D11a} can be written as
\begin{equation}
\left[\frac{T - (\nu_1 + \lambda \nu_2) M \, \ln M}{(\nu_1 + \lambda \nu_2) M} - \ln \nu_1\right] - \frac{T - T_1}{(\nu_1 + \lambda \nu_2) M} \,
\overset{D}{\longrightarrow } Y,
\qquad
M \to \infty,
\label{E12}
\end{equation}
where $Y$ is a standard Gumbel random variable. Moreover, by using \eqref{B12} and \eqref{D6a} in \eqref{E5} we get
\begin{equation}
E\left[\frac{T - T_1}{(\nu_1 + \lambda \nu_2) M}\right] = O\left(\frac{\ln M}{M^{\lambda - 1}}\right),
\qquad
M \to \infty.
\label{E13}
\end{equation}
Since $\lambda > 1$, formula \eqref{E13} implies that, as $M \to \infty$,
\begin{equation}
\frac{T - T_1}{(\nu_1 + \lambda \nu_2) M} \to 0
\qquad \text{in the $L_1$ sense.}
\label{E14}
\end{equation}
Hence the above convergence is also in probability and, consequently, in distribution. Therefore, by using \eqref{E14} in \eqref{E12} we obtain
immediately the limiting distribution of the random variable $T$ as $M \to \infty$:

\smallskip

\textbf{Theorem 8.} Let $T$ be the number of trials required to detect all Group 1 and Group 2 coupons. Then,
\begin{equation*}
\frac{T - (\nu_1 + \lambda \nu_2) M \, \ln M}{(\nu_1 + \lambda \nu_2) M} - \ln \nu_1 \,
\overset{D}{\longrightarrow } Y
\qquad \text{as }\;
M \to \infty,
\end{equation*}
where $Y$ is a standard Gumbel random variable.

\section{APPENDIX}

\subsection{An alternative derivation of formula \eqref{APP3}}
%
%
%
%
%
%

Consider the $g$-dimensional Markov chain
\begin{equation}
X(\tau) = \big(X_1(\tau), \dots, X_g(\tau)\big),
\qquad
\tau = 0, 1, 2, \dots,
\label{APPA0b}
\end{equation}
where $X_j(\tau)$, $j = 1, \dots, g$, is the number of (different) Group $j$ coupons detected after $\tau$ trials.

For typographical convenience we write
\begin{equation}
m = (m_1, m_2, \dots, m_g)
\qquad \text{and} \qquad
\label{APPA0c}
e_j = (\delta_{1j}, \delta_{2j}, \dots, \delta_{gj}),
\end{equation}
where $\delta_{ij}$ is the Kronecker delta (thus $e_j$, $j = 1, \dots, g$, is the $g$-tuple whose $j$-th entry is $1$,
while all other entries are $0$). Then the transition probabilities of $X(\tau)$ are
\begin{equation}
\left.
  \begin{array}{ccc}
    P\big\{X(\tau + 1) = m + e_j \, \big| \, X(\tau) = m \big\} & = & (M_j - m_j) p_j
    \\
    \\
    P\big\{X(\tau + 1) = m  \, \big| \, X(\tau) = m \big\}  & = & \ m_1 p_1 + \cdots + m_g p_g
    \\
  \end{array}
\right\}
\label{APPA1}
\end{equation}
for $0 \leq m_j \leq M_j$, $j = 1, \dots, g$.

We, also, introduce the quantity
\begin{equation}
u(m) := P\{T_1 = T_{\min} \, | \, X(0) = m\},
\qquad
0 \leq m_j \leq M_j, \ j = 1 \dots, g.
\label{APPA2a}
\end{equation}
We need to compute $u(0,0, \dots, 0) = P\{T_1 = T_{\min}\}$.

Using the transition probabilities given in \eqref{APPA1} we get by the law of total probability
\begin{equation*}
u(m) = \sum_{j=1}^g u(m + e_j) (M_j - m_j) p_j + u(m) (m_1 p_1  + \cdots + m_g p_g)
\end{equation*}
or
\begin{equation}
\sum_{j=1}^g p_j (M_j - m_j) \big[u(m + e_j) - u(m)\big]  = 0
\label{APPA2}
\end{equation}
for $0 \leq m_j \leq M_j$, and $j = 1, \dots, g$, which is a linear first-order (since at any time step each
$X_j(\tau)$ can only increase by $1$ or stay the same) partial difference equation for $u(m)$.

Furthermore, it is clear that $u(m)$ satisfies the boundary conditions
\begin{align}
& \ \ u(m) = 1
\qquad
\text{when }\; m_1 = M_1 \; \text{ and }\; 0 \leq m_j \leq M_j - 1 \; \text{ for }\; j \ne 1,
\nonumber
\\
& \text{while for }\; k = 2, \dots, g
\label{APPA3}
\\
& \ \ u(m) = 0
\qquad
\text{when }\; m_k = M_k \; \text{ and }\; 0 \leq m_j \leq M_j - 1 \; \text{ for }\; \ j \ne k.
\nonumber
\end{align}
It is easy to see that the resulting boundary value problem \eqref{APPA2}-\eqref{APPA3} has a unique solution $u(m)$ (incidentally,
the continuous analog of \eqref{APPA2}-\eqref{APPA3}, namely the problem
$p_1 (M_1 - x_1) \, \partial_{x_1} U(x) + \cdots + p_g (M_2 - x_g) \, \partial_{x_g} U(x) = 0$, for $x = (x_1, \dots, x_g)$ in the open box
$B = (0, M_1) \times \cdots \times (0, M_g)$,
with boundary conditions $U(x) = 1$ for $x_1 = M_1$, and $U(x) = 0$
for $x_k = M_k$, $k = 2, \dots, g$, has no solution since the hyperplanes $\{x_1 = M_1\}, \dots, \{x_g = M_g\}$ on which the Cauchy data is given,
are characteristic hypersurfaces).

If we make the simple change of variables
\begin{equation}
n = (n_1, \dots, n_g) := (M_1 - m_1, \dots, M_g - m_g)
\qquad \text{and} \qquad
\label{APPA0d}
v(n) := u(m),
\end{equation}
then the boundary value problem \eqref{APPA2}-\eqref{APPA3} can be written equivalently as
\begin{equation}
\sum_{j=1}^g p_j n_j \big[v(n - e_j) - v(n)\big]  = 0
\label{APPA2c}
\end{equation}
for $0 \leq n_j \leq M_j$, and $j = 1, \dots, g$, with the boundary conditions
\begin{align}
& \ \ v(n) = 1
\qquad
\text{when }\; n_1 = 0 \; \text{ and }\; 1 \leq n_j \leq M_j \; \text{ for }\; j \ne 1,
\nonumber
\\
& \text{while for }\; k = 2, \dots, g
\label{APPA3c}
\\
& \ \ v(n) = 0
\qquad
\text{when }\; n_k = 0 \; \text{ and }\; 1 \leq n_j \leq M_j \; \text{ for }\; \ j \ne k.
\nonumber
\end{align}
Recall that our final goal is to compute
\begin{equation}
P\{T_1 = T_{\min}\} = v(M_1, \dots, M_g).
\label{APPA4a}
\end{equation}

We will determine the solution $v(n)$ by the method of separation of variables. We first look
for solutions of \eqref{APPA2c} of the form
\begin{equation}
v(n) = \psi_1(n_1) \cdots \psi_g(n_g).
\label{APPA4}
\end{equation}
Substituting \eqref{APPA4} in equation \eqref{APPA2c} yields
\begin{equation}
\sum_{j=1}^g p_j n_j \, \frac{\psi_j(n_j - 1) - \psi_j(n_j)}{\psi_j(n_j)} = 0
\label{APPA5}
\end{equation}
and, since the first term in the sum is a function of $n_1$ alone, the second term is a function of $n_2$ alone, and so on
we must have
\begin{equation}
p_j n_j \, \frac{\psi_j(n_j - 1) - \psi_j(n_j)}{\psi_j(n_j)} = \lambda_j,
\qquad
j = 1, \dots, g,
\label{APPA6a}
\end{equation}
with
\begin{equation}
\lambda_1 + \cdots + \lambda_g = 0,
\qquad \text{i.e.} \qquad
\lambda_1 = -(\lambda_2 + \cdots + \lambda_g),
\label{APPA6b}
\end{equation}
where the constants $\lambda_j$ are the separation parameters (notice that, due to \eqref{APPA6b} there are $g-1$ independent such parameters).

Equation \eqref{APPA6a} implies
\begin{equation}
\frac{\psi_j(n_j - 1)}{\psi_j(n_j)} = \frac{p_j n_j + \lambda_j}{p_j n_j},
\label{APPA7}
\end{equation}
hence
\begin{equation}
\psi_j(n_j; \lambda_j) = p_j^{n_j} n_j!
\prod_{k = 1}^{n_j} \frac{1}{\lambda_j + p_j k},
\qquad
1 \leq n_j \leq M_j,
\label{APPA8}
\end{equation}
where without loss of generality we have assumed
\begin{equation}
\psi_j(0; \lambda_j) = 1
\qquad \text{for }\;
j = 1, \dots, g.
\label{APPA10}
\end{equation}
The quantity $\psi_j(n_j; \lambda_j)$ of \eqref{APPA8}, viewed as a rational function of $\lambda_j$, has the
partial fraction expansion
\begin{equation}
\psi_j(n_j; \lambda_j)
= \sum_{k = 1}^{n_j} (-1)^{k-1} \binom{n_j}{k} \frac{p_j k}{\lambda_j + p_j k},
\qquad
1 \leq n_j \leq M_j.
\label{APPA8a}
\end{equation}

We will now try to write the solution $v(n)$ of the boundary value problem \eqref{APPA2c}-\eqref{APPA3c} as a superposition of the special solutions
$\psi_1(n_1; \lambda_1) \cdots \psi_g(n_g; \lambda_g)$, namely we will try to express $v(n)$ as
\begin{equation}
v(n) = \int_D \psi_1(n_1; \lambda_1) \cdots \psi_g(n_g; \lambda_g) \, h(\lambda_2, \dots, \lambda_g) \, d\lambda_2 \cdots d\lambda_g
\label{APPA11}
\end{equation}
for a suitable function $h(\lambda_2, \dots, \lambda_g)$ and a suitable domain of integration $D$. Let as emphasize that $h$ and $D$ are not unique,
but one convenient pair $(h, D)$ will be enough for us.

Since \eqref{APPA2c} is a linear equation and $v(n)$ is unique, we basically need to find $h$ and $D$ such that the boundary conditions \eqref{APPA3c}
are satisfied, that is (in view of \eqref{APPA11} and \eqref{APPA10})
\begin{equation}
\int_D \psi_2(n_2; \lambda_2) \cdots \psi_g(n_g; \lambda_g)
\, h(\lambda_2, \dots, \lambda_g) \, d\lambda_2 \cdots d\lambda_g = 1
\label{APPA12a}
\end{equation}
for $n_j = 1, \dots, M_j$, $j = 2, \dots, g$, and for $k = 2, \dots, g$
\begin{align}
&\int_D \psi_1(n_1; \lambda_1) \cdots \psi_{k-1}(n_{k-1}; \lambda_{k-1}) \, \psi_{k+1}(n_{k+1}; \lambda_{k+1}) \cdots \psi_g(n_g; \lambda_g)
\, h(\lambda_2, \dots, \lambda_g) \, d\lambda_2 \cdots d\lambda_g
\nonumber
\\
& = 0
\label{APPA12b}
\end{align}
for $n_j = 1, \dots, M_j$, $j \ne k$.

From formula \eqref{APPA8} or \eqref{APPA8a} we see that the quantity $\psi_j(n_j; \lambda_j)$ is a rational
functions of $\lambda_j$ with simple poles lying
on the negative real axis. Motivated by this simple fact, and after some trial and error, we came up with the
following choice for $D$ and $h$:
\begin{equation}
D = \Gamma_2 \times \cdots \times \Gamma_g,
\label{APPA9a}
\end{equation}
where, for each $j = 2, \dots, g$ the set $\Gamma_j$ is the oriented straight line of the complex $\lambda_j$-plane (parallel to the imaginary axis)
described by the parametric equation
\begin{equation}
\lambda_j(t) = -\varepsilon + it,
\qquad
-\infty < t < \infty,
\label{APPA9b}
\end{equation}
where $\varepsilon$ is a fixed (real) number in the interval $(0, p_{\min} / g)$ with $p_{\min} = \bigwedge_{j = 2}^g p_j$;
\begin{equation}
h(\lambda_2, \dots, \lambda_g) = \frac{(-1)^{g-1}}{(2 \pi i)^{g-1}} \cdot \frac{1}{\lambda_2 \cdots \lambda_g}.
\label{APPA9c}
\end{equation}
Using the above $D$ and $h$ in \eqref{APPA11} gives
\begin{equation}
v(n) = \frac{(-1)^{g-1}}{(2 \pi i)^{g-1}} \int_{\Gamma_g} \cdots \int_{\Gamma_2}
\psi_1(n_1; \lambda_1) \psi_2(n_2; \lambda_2) \cdots \psi_g(n_g; \lambda_g) \,
\frac{d\lambda_2 \cdots d\lambda_g}{\lambda_2 \cdots \lambda_g},
\label{APPA21}
\end{equation}
where $\lambda_1$ is given by \eqref{APPA6b}.

First, let us assume that $n_j \ne 0$ for every $j = 1, \dots, g$. Then, in view of \eqref{APPA8a} we have
\begin{align}
& \psi_1(n_1; \lambda_1) \psi_2(n_2; \lambda_2) \cdots \psi_g(n_g; \lambda_g)
\nonumber
\\
& = (-1)^g \sum_{k_g = 1}^{n_g} \cdots \sum_{k_1 = 1}^{n_1} (-1)^{k_1 + \cdots + k_g}
\binom{n_1}{k_1} \cdots \binom{n_g}{k_g}
\frac{(k_1 p_1) \cdots (k_g p_g)}{(\lambda_1 + k_1 p_1) \cdots (\lambda_g + k_g p_g)}.
\label{APPA22}
\end{align}
By using \eqref{APPA22} we can see that $v(n)$ of \eqref{APPA21} is a sum of integrals of the form
\begin{equation}
-\frac{1}{(2 \pi i)^{g-1}} \int_{\Gamma_g} \cdots \int_{\Gamma_2}
\frac{(k_1 p_1) \cdots (k_g p_g)}{(\lambda_1 + k_1 p_1) (\lambda_2 + k_2 p_2) \cdots (\lambda_g + k_g p_g)}
\cdot \frac{d\lambda_2 \cdots d\lambda_g}{\lambda_2 \cdots \lambda_g}.
\label{APPA23}
\end{equation}
It is not hard to calculate the above iterated integral. We first integrate with respect to $\lambda_2$. Keeping in mind that
$\lambda_1 = -(\lambda_2 + \cdots + \lambda_g)$ and that $\Re(\lambda_j) = \varepsilon \in (0, p_{\min} / g)$ for $j = 2, \dots, g$, we see that the
integrand, viewed as a function of $\lambda_2$ has three simple poles: (i) at $\lambda_2 = -k_2 p_2$, (ii) at $\lambda_2 = 0$,
and (iii) at $\lambda_2 = k_1 p_1 - (\lambda_3 + \cdots + \lambda_g)$. Due to the choice of $\varepsilon$, the third pole has a positive real part.
Hence, the only pole which lies left of the line $\Gamma_2$ is the pole at $\lambda_2 = -k_2 p_2$. It follows that
\begin{align}
&\frac{1}{2 \pi i} \int_{\Gamma_2}
\frac{(k_1 p_1) \cdots (k_g p_g)}{(\lambda_1 + k_1 p_1) (\lambda_2 + k_2 p_2) \cdots (\lambda_g + k_g p_g)}
\cdot \frac{d\lambda_2}{\lambda_2 \cdots \lambda_g}
\nonumber
\\
= &-\frac{(k_1 p_1) (k_3 p_3) \cdots (k_g p_g)}
{(k_1 p_1 + k_2 p_2 - \lambda_3 - \cdots - \lambda_g) (\lambda_3 + k_3 p_3) \cdots (\lambda_g + k_g p_g) \lambda_3 \cdots \lambda_g},
\label{APPA24}
\end{align}
i.e. the value of the integral is the residue of the integrand at $\lambda_2 = -k_2 p_2$. Noticing that the resulting value in \eqref{APPA24}
is of the same form as the integrand with the $\lambda_2$ factors missing, makes it very easy to finish the calculation of the iterated integral
of \eqref{APPA23}. The result is
\begin{align}
&\frac{1}{(2 \pi i)^{g-1}} \int_{\Gamma_g} \cdots \int_{\Gamma_2}
\frac{(k_1 p_1) \cdots (k_g p_g)}{(\lambda_1 + k_1 p_1) (\lambda_2 + k_2 p_2) \cdots (\lambda_g + k_g p_g)}
\cdot \frac{d\lambda_2 \cdots d\lambda_g}{\lambda_2 \cdots \lambda_g}
\nonumber
\\
= & (-1)^{g-1} \frac{k_1 p_1}
{k_1 p_1 + k_2 p_2 + \cdots + k_g p_g}
\label{APPA25}
\end{align}
and by using \eqref{APPA25} and \eqref{APPA22} in \eqref{APPA21} we obtain
\begin{equation}
v(n) = (-1)^g \sum_{k_g = 1}^{n_g} \cdots \sum_{k_1 = 1}^{n_1} (-1)^{k_1 + \cdots + k_g}
\binom{n_1}{k_1} \cdots \binom{n_g}{k_g}
\frac{k_1 p_1} {k_1 p_1 + k_2 p_2 + \cdots + k_g p_g}.
\label{APPA26}
\end{equation}

Next, suppose that $n_2 = 0$. Then, since $\psi_2(0; \lambda_2) \equiv 1$ no factor of the form $\lambda_2 + k_2 p_2$ appears in the denominator
of the integrand in \eqref{APPA21}, hence when we integrate with respect to $\lambda_2$ we will get $0$. For exactly the same reason, if
$n_k = 0$ for some $k = 3, \dots, g$, then the integral in \eqref{APPA21} will vanish. Therefore, $v(n)$ of \eqref{APPA21} satisfies the boundary
conditions \eqref{APPA12b}.

Finally, if $D$ and $h$ are given by \eqref{APPA9a}-\eqref{APPA9b} and \eqref{APPA9c} respectively, then the integral in \eqref{APPA12a}
is equal to
\begin{equation}
(-1)^{g-1} \prod_{j=2}^g \frac{1}{2 \pi i} \int_{\Gamma_j} \frac{\psi_j(n_j; \lambda_j)}{\lambda_j}
\, d\lambda_j.
\label{APPA27}
\end{equation}
By formula \eqref{APPA8} we have
\begin{equation}
\frac{1}{2 \pi i} \int_{\Gamma_j} \frac{\psi_j(n_j; \lambda_j)}{\lambda_j} \, d\lambda_j
= \frac{1}{2 \pi i} \int_{\Gamma_j}
\left(\prod_{k = 1}^{n_j} \frac{k p_j}{\lambda_j + k p_j}\right) \frac{d\lambda_j}{\lambda_j} = -1,
\label{APPA28}
\end{equation}
where the second equality in \eqref{APPA28} follows immediately by the observation that the only pole of the integrand on the right of
$\Gamma_j$ is located at $\lambda_j = 0$ (and the residue there is $1$, however due to the orientation of $\Gamma_j$ the value of the integral
is $-1$). Formula \eqref{APPA28} implies that the value of the quantity in \eqref{APPA27} is $1$ and this verifies that $v(n)$ of \eqref{APPA21}
satisfies the boundary conditions \eqref{APPA12a}. We have thus shown that $v(n)$ of \eqref{APPA26} is the solution of the boundary value problem
\eqref{APPA2c}-\eqref{APPA3c}. Therefore the proof of formula \eqref{APP3} is completed by setting $n = (M_1, \dots, M_g)$ in \eqref{APPA26} and
invoking \eqref{APPA4a}.
\hfill $\blacksquare$

\subsection{Proof of formula \eqref{AUX30}}
Let $\alpha \in (0, 1)$. For typographical convenience we set
\begin{equation}
U(N; \alpha) := e^{\ln^{\alpha} N}
\label{AUX8}
\end{equation}
(so that for any constant $\beta > 0$ we have $\ln^{\beta}N << U(N; \alpha) << N^{\beta}$ as $N \to \infty$)
and then we write \eqref{AUX7} as
\begin{equation}
E\left[S_N^{(r)}\right] = N^r \ln^r N \left[I_1(N) + I_2(N)\right],
\label{AUX9}
\end{equation}
where
\begin{equation}
I_1(N) := \int_0^{U(N; \alpha)} \left(1 - \frac{x}{N}\right)^{N-1} \left(1 - \frac{\ln x}{\ln N}\right)^r dx
\label{AUX10a}
\end{equation}
and
\begin{equation}
I_2(N) := \int_{U(N; \alpha)}^N \left(1 - \frac{x}{N}\right)^{N-1} \left(1 - \frac{\ln x}{\ln N}\right)^r dx.
\label{AUX10b}
\end{equation}
We will first estimate $I_2(N)$ as $N \to \infty$.
\begin{align}
0 < I_2(N) &< \left(1 - \frac{U(N; \alpha)}{N}\right)^{N-1}
\int_{U(N; \alpha)}^N \left(1 - \frac{\ln x}{\ln N}\right)^r dx
\nonumber
\\
&< N \left[\left(1 - \frac{U(N; \alpha)}{N}\right)^{\frac{N-1}{U(N; \alpha)}}\right]^{U(N; \alpha)}
\sim N e^{-U(N; \alpha)}.
\label{AUX11}
\end{align}
In particular \eqref{AUX11} implies
\begin{equation}
I_2(N) = o\left(\frac{1}{N^{\kappa}}\right)
\qquad \text{for any }\;
\kappa > 0.
\label{AUX12}
\end{equation}
To estimate $I_1(N)$ we first notice that for $0 \leq x \leq U(N; \alpha) = e^{\ln^{\alpha} N}$ we have
\begin{equation}
(N-1) \ln\left(1 - \frac{x}{N}\right) = -(N-1) \left[\frac{x}{N} + O\left(\frac{x^2}{N^2}\right)\right]
= -x + o\left(\frac{1}{N^{\theta}}\right)
\quad \text{for any }\;
\theta \in (0, 1).
\label{AUX13}
\end{equation}
Exponentiating \eqref{AUX13} yields
\begin{equation}
\left(1 - \frac{x}{N}\right)^{N-1} = e^{-x} \left[1 + o\left(\frac{1}{N^{\theta}}\right)\right]
\quad \text{for any }\;
\theta \in (0, 1).
\label{AUX14}
\end{equation}
We, then, substitute \eqref{AUX14} in \eqref{AUX10a} and obtain
\begin{equation}
I_1(N) = \int_0^{U(N; \alpha)} e^{-x} \left(1 - \frac{\ln x}{\ln N}\right)^r dx \,
+ \, o\left(\frac{1}{N^{\theta}}\right) \int_0^{U(N; \alpha)} e^{-x} \left(1 - \frac{\ln x}{\ln N}\right)^r dx.
\label{AUX15}
\end{equation}
Now,
\begin{equation}
\int_0^{U(N; \alpha)} e^{-x} \left(1 - \frac{\ln x}{\ln N}\right)^r dx
\leq \int_0^1 e^{-x} \left(1 - \frac{\ln x}{\ln N}\right)^r dx
+ \int_1^{U(N; \alpha)} e^{-x} dx = O(1)
\label{AUX16}
\end{equation}
as $N \to \infty$. Thus, by using \eqref{AUX16} in \eqref{AUX15} we arrive at
\begin{equation}
I_1(N) = \int_0^{U(N; \alpha)} e^{-x} \left(1 - \frac{\ln x}{\ln N}\right)^r dx \,
+ \, o\left(\frac{1}{N^{\theta}}\right)
\qquad \text{for any }\;
\theta \in (0, 1).
\label{AUX17}
\end{equation}
Since $\ln x \to -\infty$ as $x \to 0^+$, we need to estimate the ``bottom tail" of the integral in
\eqref{AUX17}. We have
\begin{equation}
0 < \int_0^{U(N; \alpha)^{-1}} e^{-x} \left(1 - \frac{\ln x}{\ln N}\right)^r dx
< \left(\frac{2^r}{\ln^r N} + \frac{2^r}{e^{\ln^{\alpha} N}}\right)\int_0^{U(N; \alpha)^{-1}} (-\ln x)^r dx
\label{AUX18}
\end{equation}
(recall that $U(N; \alpha)^{-1} = e^{-\ln^{\alpha} N}$, where $0 < \alpha < 1$). Integration by parts gives
\begin{equation}
\int_0^{U(N; \alpha)^{-1}} (-\ln x)^r dx = e^{-\ln^{\alpha} N} (\ln N)^{r \alpha} \, [1 + o(1)],
\label{AUX19}
\end{equation}
hence, by using \eqref{AUX19} in \eqref{AUX18} we obtain
\begin{equation}
\int_0^{U(N; \alpha)^{-1}} e^{-x} \left(1 - \frac{\ln x}{\ln N}\right)^r dx
= o\left(e^{-\ln^{\alpha} N}\right).
\label{AUX20}
\end{equation}
In view of \eqref{AUX20}, formula \eqref{AUX17} implies
\begin{equation}
I_1(N) = \int_{U(N; \alpha)^{-1}}^{U(N; \alpha)} e^{-x} \left(1 - \frac{\ln x}{\ln N}\right)^r dx \,
+ \, o\left(e^{-\ln^{\alpha} N}\right).
\label{AUX21}
\end{equation}
Finally, by substituting \eqref{AUX12} and \eqref{AUX21} in \eqref{AUX9} we get
\begin{equation}
E\left[S_N^{(r)}\right] = N^r (\ln N)^r J(N; \alpha),
\label{AUX22a}
\end{equation}
where
\begin{equation}
J(N; \alpha) :=
\int_{U(N; \alpha)^{-1}}^{U(N; \alpha)} e^{-x} \left(1 - \frac{\ln x}{\ln N}\right)^r dx \,
+ \, o\left(e^{-\ln^{\alpha} N}\right)
\label{AUX22b}
\end{equation}
as $N \to \infty$.

Before we continue let us recall that, for any given $r > 0$, the function $h(y) := (1 - y)^r$ has the Taylor series expansion
\begin{equation}
h(y) = (1 - y)^r = \sum_{k=0}^{\infty} (-1)^k \binom{r}{k} y^k,
\qquad
y \in [-1, 1],
\label{AUX23a}
\end{equation}
where
\begin{equation}
\binom{r}{0} = 1
\qquad \text{and} \qquad
\binom{r}{k} = \frac{r (r-1) \cdots (r-k+1)}{k!}
\quad \text{for }\;
k = 1, 2, \dots \;.
\label{AUX23b}
\end{equation}
Thus, if $n \geq 1$ is a fixed integer, then
\begin{equation}
h(y) = (1 - y)^r = \sum_{k=0}^n (-1)^k \binom{r}{k} y^k + \frac{h^{(n+1)}(\xi)}{(n+1)!} y^{n+1},
\label{AUX24}
\end{equation}
where $\xi$ lies between $0$ and $y$.

Going back to formula \eqref{AUX22b} we see that due to the limits of integration the dummy variable $x$ satisfies
\begin{equation}
-\frac{\ln^{\alpha} N}{\ln N} \leq \frac{\ln x}{\ln N} \leq \frac{\ln^{\alpha} N}{\ln N}.
\label{AUX25}
\end{equation}
Hence, if we set $y = \frac{\ln x}{\ln N}$ in \eqref{AUX24}, then the quantity $\frac{h^{(n+1)}(\xi)}{(n+1)!}$
is bounded and we get
\begin{equation}
\left(1 - \frac{\ln x}{\ln N}\right)^r = \sum_{k=0}^n (-1)^k \binom{r}{k} \left(\frac{\ln x}{\ln N}\right)^k
+ O\left(\frac{\ln^{\alpha(n+1)} N}{\ln^{n+1} N}\right)
\label{AUX26}
\end{equation}
uniformly in $x$, as long as the range of values of $x$ is given by \eqref{AUX25}.

It follows from \eqref{AUX26} that if (for our given $n$) we choose an $\alpha$ so that
\begin{equation}
0 < \alpha < \frac{1}{n+1},
\label{AUX26a}
\end{equation}
then
\begin{equation}
\left(1 - \frac{\ln x}{\ln N}\right)^r = \sum_{k=0}^n (-1)^k \binom{r}{k} \left(\frac{\ln x}{\ln N}\right)^k
+ o\left(\frac{1}{\ln^n N}\right)
\label{AUX27}
\end{equation}
again uniformly in $x$, within the range of values given by \eqref{AUX25}. Thus, we can substitute \eqref{AUX27} in \eqref{AUX22b} and get
\begin{equation}
J(N; \alpha) = \sum_{k=0}^n \binom{r}{k} \frac{(-1)^k}{\ln^k N}
\int_{U(N; \alpha)^{-1}}^{U(N; \alpha)} e^{-x} (\ln x)^k dx
+ o\left(\frac{1}{\ln^n N}\right)
\label{AUX28}
\end{equation}
as $N \to \infty$.

Next, we observe that in the same way we derived \eqref{AUX20} we can also get
\begin{equation}
\int_0^{U(N; \alpha)^{-1}} e^{-x} (\ln x)^k dx
= o\left(e^{-\ln^{\alpha} N} \ln^k N\right).
\label{AUX29a}
\end{equation}
Also, it is easy to see that
\begin{equation}
\int_{U(N; \alpha)}^{\infty} e^{-x} (\ln x)^k dx
= O\left(e^{-U(N; \alpha)} (\ln N)^{\alpha k}\right)
\label{AUX29b}
\end{equation}
(recall that $U(N; \alpha) = e^{\ln^{\alpha} N}$).
Therefore, by using \eqref{AUX29a} and \eqref{AUX29b} in \eqref{AUX28} we obtain
\begin{equation}
J(N; \alpha) = \sum_{k=0}^n \binom{r}{k} \frac{(-1)^k}{\ln^k N}
\int_0^{\infty} e^{-x} (\ln x)^k dx
+ o\left(\frac{1}{\ln^n N}\right),
\qquad
N \to \infty,
\label{AUX30a}
\end{equation}
and, finally, by substituting \eqref{AUX30a} in \eqref{AUX22a} we arrive at \eqref{AUX30}.
\hfill $\blacksquare$

\subsection{Proof of formula \eqref{A6f}}
We start by introducing the generating functions
\begin{equation}
G(z) := E\left[z^{-S(\theta)}\right]
= 1 - (z-1) \int_0^{\infty} e^{-(z-1) t} \left[1 - \left(1 - e^{-\theta t} \right) \left(1 - e^{-(1 - \theta) t / N} \right)^N\right] dt.
\label{APPD1}
\end{equation}
%
Notice that if
\begin{equation}
\Re\{z\} > 1 - \frac{1 - \theta}{N},
\label{APPD00}
\end{equation}
the integral appearing in \eqref{APPD1} is absolutely convergent.

We will derive formula \eqref{A6f} via characteristic functions. Let us fix a $\xi \in \mathbb{R}$ and set
\begin{equation}
\zeta := e^{-i \xi}.
\label{APPD2}
\end{equation}
Then, in view of \eqref{APPD1}, the characteristic function of $[(1 - \theta) S(\theta) - N \ln N] /N$ is
\begin{equation}
\phi_N(\xi) = E\left[\zeta^{\,-\frac{(1 - \theta) S(\theta) - N \ln N}{N}}\right]
= \zeta^{\,\ln N} E\left[\left(\zeta^{\frac{(1 - \theta)}{N}}\right)^{-S(\theta)}\right]
= \zeta^{\,\ln N} G\left(\zeta^{\frac{(1 - \theta)}{N}}\right).
\label{APPD3}
\end{equation}
Now,
\begin{equation}
\zeta^{\frac{(1 - \theta)}{N}} = e^{-\frac{i(1 - \theta) \xi}{N}} = 1 - \frac{i(1 - \theta) \xi}{N} + O\left(\frac{1}{N^2}\right),
\qquad
N \to \infty.
\label{APPD4}
\end{equation}
In particular $z = \zeta^{(1 - \theta) / N}$ satisfies \eqref{APPD00} for all $N$ sufficiently large.

Next, by using \eqref{APPD1} and \eqref{APPD4} in \eqref{APPD3} we get
\begin{equation}
\zeta^{\,-\ln N} \phi_N(\xi)
= 1 + \left[\frac{i(1 - \theta) \xi}{N} + O\left(\frac{1}{N^2}\right)\right] \big[\chi_1(N) + \chi_2(N)\big],
\label{APPD5}
\end{equation}
where
\begin{equation}
\chi_1(N) :=
\int_0^{\infty} e^{-\left(\zeta^{\frac{(1 - \theta)}{N}} - 1 \right) t} \left[1 - \left(1 - e^{-(1 - \theta) t / N} \right)^N\right] dt
\label{APPD6a}
\end{equation}
and
\begin{equation}
\chi_2(N) :=
\int_0^{\infty} e^{-\left(\zeta^{\frac{(1 - \theta)}{N}} - 1 \right) t} e^{-\theta t} \left(1 - e^{-(1 - \theta) t / N} \right)^N dt.
\label{APPD6b}
\end{equation}
Regarding $\chi_2(N)$, in the same way we got formula \eqref{A3i} from \eqref{A3h} we can obtain
\begin{equation}
\chi_2(N) = O\left(e^{-\varepsilon N}\right),
\qquad
N \to \infty,
\label{APPD7}
\end{equation}
for any sufficiently small $\varepsilon > 0$.

Now, using \eqref{APPD4} in \eqref{APPD6a} yields
\begin{equation}
\chi_1(N)
= \int_0^{\infty} e^{[i(1 - \theta) \xi \, + \, O(N^{-1})] \, t / N}
\left[1 - \left(1 - e^{-(1 - \theta) t / N} \right)^N\right] dt
\label{APPD8}
\end{equation}
or, after the substitution $s = (1 - \theta) t / N$ in the above integral
\begin{equation}
\chi_1(N)
= \frac{N}{1 - \theta} \int_0^{\infty} e^{[i \xi \, + \, O(N^{-1})] \, s}
\left[1 - \left(1 - e^{- s} \right)^N\right] ds.
\label{APPD9}
\end{equation}
Therefore, by substituting \eqref{APPD7} and \eqref{APPD9} in \eqref{APPD5} we obtain
\begin{align}
\zeta^{\,-\ln N} \phi_N(\xi)
&= 1 + \left[i \xi + O\left(\frac{1}{N}\right)\right]
\int_0^{\infty} e^{[i \xi \, + \, O(N^{-1})] \, s}
\left[1 - \left(1 - e^{- s} \right)^N\right] ds
\nonumber
\\
&= 1 + i \xi \int_0^{\infty} e^{i \xi s}
\left[1 - \left(1 - e^{- s} \right)^N\right] ds + O\left(\frac{1}{N}\right)
\label{APPD10}
\end{align}
as $N \to \infty$.

Let $A_N$ be a (real) quantity which grows to $\infty$ with $N$ so that
\begin{equation}
\frac{N}{e^{A_N}} = o\left(\frac{1}{N}\right)
\qquad \text{as }\;
N \to \infty
\label{APPD11}
\end{equation}
(we do not need to be more specific about $A_N$). Noticing that by \eqref{APPD2} we have $\zeta^{\,-\ln N} = N^{i \xi}$,
we rewrite \eqref{APPD10} as
\begin{equation}
N^{i \xi} \phi_N(\xi) = 1 + K_1(N) + K_2(N) + O\left(\frac{1}{N}\right),
\label{APPD10a}
\end{equation}
where
\begin{equation}
K_1(N) := \int_0^{A_N} i \xi \, e^{i \xi s} \left[1 - \left(1 - e^{- s} \right)^N\right] ds
\label{APPD10b}
\end{equation}
and
\begin{equation}
K_2(N) := \int_{A_N}^{\infty} i \xi \, e^{i \xi s} \left[1 - \left(1 - e^{- s} \right)^N\right] ds.
\label{APPD10c}
\end{equation}
Applying integration by parts in \eqref{APPD10b} yields
\begin{equation}
K_1(N) = e^{i \xi A_N} \left[1 + \left(1 - e^{-A_N}\right)^N\right] - 1
+ N \int_0^{A_N} e^{i \xi s} \left(1 - e^{- s} \right)^{N-1} e^{- s} ds,
\label{APPD12}
\end{equation}
which, in view of \eqref{APPD11}, implies
\begin{equation}
K_1(N) = e^{i \xi A_N} - 1
+ N\int_0^{A_N} e^{i \xi s} \left(1 - e^{- s} \right)^{N-1} e^{- s} ds
+ o\left(\frac{1}{N}\right).
\label{APPD13}
\end{equation}
Next, by substituting $s = \ln N - \ln x$ in the integral appearing in the right-hand side of \eqref{APPD13} we obtain
\begin{equation}
K_1(N) = e^{i \xi A_N} - 1
+ N^{i \xi} \int_{N e^{-A_N}}^N x^{-i \xi} \left(1 - \frac{x}{N} \right)^{N-1} dx
+ o\left(\frac{1}{N}\right).
\label{APPD14}
\end{equation}
We, then, use \eqref{APPD14} in \eqref{APPD10a} and get
\begin{equation}
N^{i \xi} \phi_N(\xi) = e^{i A_N \xi}
+ N^{i \xi} \int_{N e^{-A_N}}^N x^{-i \xi} \left(1 - \frac{x}{N} \right)^{N-1} dx + K_2(N) + O\left(\frac{1}{N}\right).
\label{APPD15}
\end{equation}

Let us, now, turn our attention to the integral $K_2(N)$ of formula \eqref{APPD10c}.

Assume first that $\xi > 0$. We complexify the dummy variable $s$ by setting
$z = s + i\tau$ and for $N$ (temporarily) fixed we choose $R > A_N$ and consider the close contour $C_R$ formed by (i) the interval $[A_N, R]$ of the
(real) $s$-axis, (ii) the circular arc $R e^{i \theta}$, $0 \leq \theta \leq \arccos(A_N / R)$, and (iii) the line segment $A_N + i\tau$,
$0 \leq \tau \leq \sqrt{R^2 - A_N^2}$
. Then, Cauchy's Theorem implies

\begin{equation}
\oint_{C_R} i \xi \, e^{i \xi z} \left[1 - \left(1 - e^{- z} \right)^N\right] dz = 0
\qquad \text{for every }\;
R > A_N.
\label{APPD16}
\end{equation}
Next (keeping $N$ fixed), we take limits in \eqref{APPD16} as $R \to \infty$. It is a standard exercise in contour integration to show that the
integral on the circular piece of $C_R$, namely on the arc $R e^{i \theta}$, $0 \leq \theta \leq \arccos(A_N / R)$, vanishes. Hence, in
view of \eqref{APPD10c}, formula \eqref{APPD16} implies
\begin{equation}
K_2(N) = - e^{i \xi A_N} \int_0^{\infty} \xi \, e^{-\xi \tau} \left[1 - \left(1 - e^{-A_N} e^{-i \tau} \right)^N\right] d\tau.
\label{APPD17}
\end{equation}
Now, we can allow $N$ to grow large. Thus, in view of \eqref{APPD11}, formula \eqref{APPD17} yields
\begin{equation}
K_2(N) = -e^{i \xi A_N} + o\left(\frac{1}{N}\right)
\label{APPD18}
\end{equation}
and, hence, by substituting \eqref{APPD18} in \eqref{APPD15} we obtain
\begin{equation}
\phi_N(\xi) = \int_{N e^{-A_N}}^N x^{-i \xi} \left(1 - \frac{x}{N} \right)^{N-1} dx + O\left(\frac{1}{N}\right).
\label{APPD19}
\end{equation}
Formula \eqref{APPD19} was obtained under the assumption that $\xi > 0$. However, if $\xi < 0$, then the same approach works if we choose the contour
$C_R$ to be the symmetric of the previous one with respect to the (real) $s$-axis. Therefore, formula \eqref{APPD19} is valid for all
$\xi \in \mathbb{R} \setminus \{0\}$, while for $\xi = 0$ formulas \eqref{APPD2} and \eqref{APPD3} imply immediately that
\begin{equation}
\phi_N(0) = 1
\qquad \text{for all }\;
N.
\label{APPD20}
\end{equation}
Finally, as in the previous subsection (see, e.g., \eqref{AUX10a} and \eqref{AUX14}), formulas \eqref{APPD19} and \eqref{APPD20} imply
\begin{equation}
\lim_{N \to \infty} \phi_N(\xi) = \int_0^{\infty} x^{-i \xi} e^{-x} dx = \Gamma(1 - i \xi)
\qquad \text{pointwise for }\; \xi \in \mathbb{R},
\label{APPD21}
\end{equation}
where $\Gamma(1 - i\xi)$ is recognized as the characteristic function of the standard Gumbel distribution.
\hfill $\blacksquare$

\end{document}